\documentclass[12pt]{amsart}
\usepackage{graphicx}
\usepackage{color}
\usepackage{setspace}
\usepackage{eufrak}
\usepackage{float}
\usepackage{etex}
\usepackage[all]{xy}
\usepackage{amssymb, amsmath}
\usepackage{eufrak}
\usepackage{ascmac}

\usepackage{mathrsfs}
\def\qed{\hfill $\Box$}
\def\proof{\noindent {\sl Proof} :\;  }
\def\t{\noindent}

\newcommand{\A}{\mathcal{A}}

\newcommand{\R}{\mathbb{R}}

\newcommand\rank{\mbox{\rm rank}\,}

\newcommand\Imag{\mbox{\rm Im} }

\def\qed{\hfill $\Box$}
\def\proof{\noindent {\sl Proof} :\;  }
\def\t{\noindent}
\def\rd{\partial}

\newcommand{\bH}{\mathbb{H}}
\newcommand{\bD}{\mathbb{D}}
\newcommand{\cbU}{\check{\mathbb{U}}}

\def\bx{\mbox{\boldmath $x$}}

\def\bu{\mbox{\boldmath $u$}}
\def\bv{\mbox{\boldmath $v$}}
\def\bw{\mbox{\boldmath $w$}}
\def\ba{\mbox{\boldmath $a$}}

\def\be{\mbox{\boldmath $e$}}

\def\bn{\mbox{\boldmath $n$}}

\def\br{\mbox{\boldmath $r$}}
\def\b0{\mbox{\boldmath $0$}}

\def\bt{\mbox{\boldmath $t$}}
\def\cbl{\check{\bv}}
\def\cbn{\check{\bn}}
\def\cbt{\check{\bt}}
\def\ckappa{\check{\kappa}}
\def\ctau{\check{\tau}}

\newtheorem{thm}{\bf Theorem}[section]

\newtheorem{lem}[thm]{\bf Lemma}
\newtheorem{prop}[thm]{\bf Proposition}

\newtheorem{rem}[thm]{\bf Remark}

\begin{document}
\title[Dual quaternions and singularities of ruled surfaces]
{Geometric algebra and singularities of ruled and developable surfaces}
%
\author[J.~Tanaka]{Junki Tanaka}
\address[J.~Tanaka]{Kobo Co, LTD, Kobe, Japan}
\email{junki.tanaka@khobho.co.jp
}
\author[T.~Ohmoto]{Toru Ohmoto}
\address[T.~Ohmoto]{Department of Mathematics,
Faculty of Science,  Hokkaido University,
Sapporo 060-0810, Japan}
\email{ohmoto@math.sci.hokudai.ac.jp}
\subjclass[2010]{58K05, 32S15, 34A09, 53A20}
\keywords{Geometric algebra, ruled surfaces, developable surfaces, singularities of smooth maps. }
\dedicatory{Dedicated to Professor Goo Ishikawa on the occasion of his 60th birthday.}
%

%
\begin{abstract} 
Any ruled surface in $\R^3$ is described as a curve of unit dual vectors 
in the algebra of dual quaternions (=the even Clifford algebra $ C\ell^+(0,3,1)$). 
Combining this classical framework and  
$\A$-classification theory of $C^\infty$ map-germs $(\R^2,0) \to (\R^3,0)$,  
we characterize local diffeomorphic types of singular ruled surfaces 
in terms of geometric invariants. 
In particular, using a theorem of G. Ishikawa, we show that 
local topological type of singular developable surfaces is completely determined 
by vanishing order of  the {\it dual torsion} $\check{\tau}$, 
that generalizes an old result of D. Mond for tangent developables of non-singular space curves. 
This work suggests that 
Geometric Algebra would be useful for studying singularities of geometric objects 
in classical Klein geometries. 
\end{abstract}

\maketitle

\setlength{\baselineskip}{15pt}

\section{Introduction}
A {\it ruled surface} in Euclidean space $\R^3$ is a surface formed 
by a $1$-parameter family of straight lines, called  {\it rulings}; 
at least partly, 
it admits a parametrization of the form $F(s, t)=\br(s)+t\be(s)$ with $|\be(s)|=1$, 
$s \in I$, $t \in \R$, where $I$ is an open interval. 
A {\it developable surface} is a ruled surface which is locally planar 
(i.e. the Gaussian curvature is constant zero). 
The parametrization map $F: I \times \R \to \R^3$ may be singular at some point $(s_0, t_0)$, 
that is, the differential $dF(s_0, t_0)$ may have rank one, 
and then the surface (= the image of $F$) has a particularly singular shape around that point. 
In this paper, we study local diffeomorphic types of of the singular surface and its bifurcations (see Fig.1).  
All maps and manifolds are assumed to be class $C^\infty$ throughout. 

\begin{figure}[h]
 \includegraphics[width=7cm]{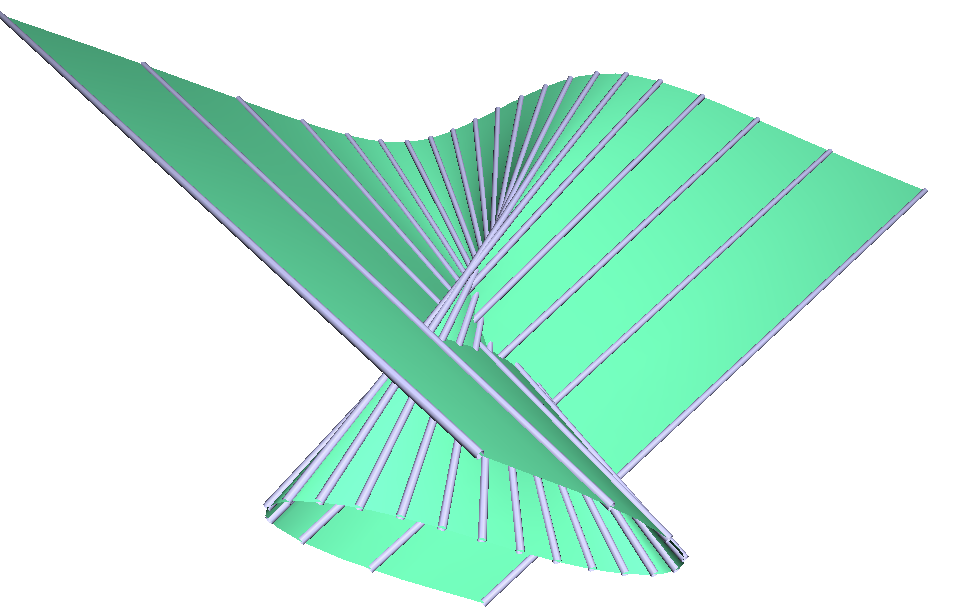}\\
 \caption{\small  
Deforming Mond's $H_2$-singularity via a family of ruled surfaces: 
the surface has two crosscaps and one triple point. 
}\label{fig1}
\end{figure}

The main feature of this paper is to combine 
classical line geometry using {\it dual quaternions}  \cite{Guggenheimmer, Hlavaty, PW, Selig}  
and $\A$-classification theory of singularities of (frontal) maps $\R^2 \to \R^3$ 
\cite{Mond85,  Ishikawa, IST, IS}. 
Here $\A$ denotes a natural equivalence relation in singularity theory of $C^\infty$ maps; 
two map-germs $f, g: (\R^2, 0) \to (\R^3,0)$ 
are {\it $\A$-equivalent} if 
there exist diffeomorphism-germs $\sigma: (\R^2,0) \to  (\R^2,0)$ 
and $\varphi:  (\R^3,0) \to  (\R^3,0)$ such that $g=\varphi\circ f \circ \sigma^{-1}$.  
We simply say the {\it $\A$-type} of a map-germ to mean its $\A$-equivalence class. 
As a weaker notion, 
 {\it topological $\A$-equivalence} is defined 
by taking $\sigma$ and $\varphi$ to be homeomorphism-germs. 
We also use the $\A$-equivalence with 
the target changes being rotations $\varphi \in SO(3)$, which 
is called  {\it rigid equivalence}  throughout the present paper. 
Our aim is to classify germs of parametrization maps $F$ of ruled surfaces in $\R^3$ 
up to $\A$-equivalence and rigid equivalence.

\subsection{Ruled surfaces} 
Geometric Algebra is a neat tool for studying motions in classical geometry;  
in the case of Euclidean $3$-space, it is the algebra of dual quaternions (e.g. Selig \cite{Selig}). 
As an application, any ruled surface in $\R^3$ is described as a curve of {\it unit dual vectors} 
$$\cbl: I \to \cbU \subset 
\bD^3, \quad 
\cbl(s)=\bv_0(s)+\varepsilon \bv_1(s) $$
Here, $\bD =\R \oplus \varepsilon \R$ with $\varepsilon^2=0$ is 
the $\R$-algebra of {\it dual numbers}, 
and $\bD^3=\R^3 \oplus \varepsilon \R^3$ is the space of {\it dual vectors}, 
and especially, 
the space of unit dual vectors is defined by 
$$\cbU:=\{\check{\bv}=\bv_0+\varepsilon \bv_1 \in \bD^3, \; 
|\bv_0|=1, \; \bv_0\cdot \bv_1=0\}$$
which is a $4$-dimensional submanifold in the $6$-dimensional space $\bD^3$. 
Obviously, $\cbU$ is diffeomorphic to the total space of the (co)tangent bundle $TS^2$, 
thus it is naturally identified with the space of oriented lines in $\R^3$, 
so $\check{\bv}$ is regarded as a $1$-parameter family of oriented lines; 
the ruled surface is parameterized by $F(s,t)=\bv_0(s)\times \bv_1(s)+t\bv_0(s)$. 
See \S 2.1 for the detail.  
In our context, as the space of ruled surfaces in $\R^3$, 
we consider the space $C^\infty(I, \cbU)$ of curves of unit dual vectors  
endowed with the Whintey $C^\infty$-topology.  

Assume that our ruled surface is {\it non-cylindrical}, 
i.e., $\bv_0'(s)\not=0$ for any $s\in I$, 
then the curve $\cbl$ admits the Frenet formula in $\bD^3$ with 
{\it complete} differential invariants, 
the {\it dual curvature} 
and the {\it dual torsion} 
$$\ckappa(s)=\kappa_0(s)+\varepsilon \kappa_1(s), 
\quad \ctau(s)=\tau_0(s)+\varepsilon \tau_1(s) \;\;\; \in \bD.$$ 
Here we may take $s$ to be the arclength of the spherical curve $\bv_0(s)$, 
that is equivalent to $\kappa_0(s)\equiv 1$, 
thus three real functions $\kappa_1, \tau_0, \tau_1$ are essential. 
In particular, $\kappa_1(s_0)=0$ if and only if 
$F$ is singular at $(s_0, t_0)$ for some (unique) $t_0$  (Lemma \ref{striction_curve}). 

We determine which $\A$-types of singular germs $(\R^2, 0) \to (\R^3,0)$ appear in 
generic families of  ruled surfaces. 
Assume that $F$ is singular at $(s_0, t_0)=(0,0)$ and $F(0,0)=0$. 
From the {\it dual Bouquet formula} of $\cbl$ at $s=0$ in $\bD^3$,  
we derive a {\it canonical} Taylor expansion of parameterization map $F$ (\S 3.2),  
where $o(n)$ denotes Landau's notation of function-germs of order greater than $n$:  
$$
\left\{
\begin{array}{ccl}
x&=& 
t-\frac{1}{2}ts^2+\frac{\tau_1(0)}{2} s^3 + o(3)\\ 
y&=&
ts-\frac{\tau_1(0)}{2}s^2-\frac{2\tau_0(0)\kappa_1'(0)+\tau_1'(0)}{6}s^3 + o(3)\\
z&=&
\frac{\kappa_1'(0)}{2}s^2+\frac{\tau_0(0)}{2}t s^2
+\frac{\kappa_1''(0)-2\tau_0(0)\tau_1(0)}{6}s^3+ o(3)
\end{array}
\right.
$$
Then we apply to the jet of $F$ the criteria 
for detecting $\A$-types of map-germs  in Mond \cite{Mond82, Mond85}.

\begin{thm}\label{thm1}
The $\A$-classification of singularities of $F$  
arising in generic at most $3$-parameter families of 
non-cylindrical ruled surfaces is given as in Table \ref{table1}; 
In particular,  for each $\A$-type in that table, 
the canonical expansion with the described condition  is regarded as  
a normal form of the jet of ruled surface-germ under rigid equivalence. 
\end{thm}

\begin{table}
{ 
$$
\begin{array}{l | l | l | l  }
&\mbox{\footnotesize normal form} & \ell &  \mbox{\footnotesize cond. at $s=s_0$}\\
\hline \hline
S_0 & (x, y^2, xy) &2& \kappa_1=0, \;\; \kappa_1'\not=0 \\
\hline
S_1^\pm& (x, y^2, y^3\pm x^2y) 
&3&  \kappa_1=\kappa_1'=0 , \;\; \tau_1\not=0, \;\; 
\kappa_1''(\kappa_1''-2\tau_0\tau_1)\gtrless 0\\
\hline
S_2& (x,y^2,y^3+ x^3y)
&4&\kappa_1=\kappa_1'=\kappa_1''=0,\; 
\kappa_1^{(3)}\tau_0\tau_1\neq 0\\ 
B_2^\pm& (x,y^2,x^2y\pm y^5)
&& \kappa_1=\kappa_1'=0,\;
\kappa_1''=2\tau_0\tau_1\neq 0,\; b_2\gtrless 0\\ 
H_2& (x,xy+y^5,y^3)
&& \kappa_1=\kappa_1'=\tau_1=0,\; \kappa_1''\neq 0,\; h_2\neq 0\\ 
\hline
S_3^\pm &(x,y^2,y^3\pm x^4y)
&5& \kappa_1=\kappa_1'=\kappa_1''=\kappa_1^{(3)}=0,\; \kappa_1^{(4)}\tau_0\tau_1 \gtrless 0\\
C_3^\pm &(x,y^2,xy^3\pm x^3y)
&& \kappa_1=\kappa_1'=\kappa_1''=\tau_0=0, \; \tau_1\neq 0, \; \kappa_1^{(3)}(\kappa_1^{(3)}-2\tau_0'\tau_1)\gtrless 0 \\
B_3^\pm &(x,y^2,x^2y\pm y^7)
&& \kappa_1=\kappa_1'=0,\;
\kappa_1''=2\tau_0\tau_1\neq 0,\; b_2=0,\ b_3\gtrless 0
\\
H_3 &(x,xy+y^7,y^3)
&& \kappa_1=\kappa_1'=\tau_1=0,\;
\kappa_1''\neq 0,\; h_2=0,\; h_3\neq 0
\\
P_3 &(x,xy+y^3,xy^2+p_4 y^4)
&& 
\kappa_1=\kappa_1'=\kappa_1'' =\tau_1= 0,\; 
\tau_0\tau_1'\neq 0, \; p_4 \not=0, 1,  \frac{1}{2},  \frac{3}{2}. 
\end{array}
$$
}
\caption{$\A$-types of singularities of ruled surfaces. 
In conditions, $\kappa_1', \kappa_1'', \cdots$ are derivatives at $s=s_0$ for short, 
e.g. $\kappa_1'$ means $\frac{d}{ds}\kappa_1(s_0)$, 
and $b_2, b_3, h_2, h_3, p_4$  are 
some polynomials of those derivatives (see \S 3.2). 
The letters  $\lessgtr, \gtrless, \pm$ are in the same order. 
In the second column, $\ell$ means $\A$-codimension of the map-germ. }
\label{table1}
\end{table}

Precisely saying, via a variant of Thom's transversality theorem (\S 3.3), 
we show that 
there exists a dense subset $\mathcal{O}$ in 
the mapping space $\mathcal{R}_W$ consisting of families of non-cylindrical  $\cbl: I\times W \to \cbU$ 
with parameter space $W$ of dimension $\le 3$ 
so that for any family belonging to $\mathcal{O}$ and for any $\lambda \in W$, 
the germ of the corresponding paramatrization map 
$F(-, \lambda): I \times \R \to \R^3$ at any point $(s_0, t_0)$ is $\A$-equivalent to either of an immersion-germ or one of singular germs in Table \ref{table1}. 

Obviously, normal forms under rigid equivalence have functional moduli: 
those are nothing but $\kappa_1(s)$, $\tau_0(s)$ and $\tau_1(s)$ having the prescribed condition 
on derivatives at $s=s_0$. 

\begin{rem}\label{rem_realization}\upshape 
({\bf Realization}) 
Izumiya-Takeuchi \cite{IT} firstly proved 
in a rigorous way that a generic singularity of ruled surfaces is only of type crosscap $S_0$, 
and  Martins and Nu\~no-Ballesteros \cite{MN} showed that 
any $\A$-simple map-germ $(\R^2, 0) \to (\R^3, 0)$ is $\A$-equivalent to a germ of ruled surface. 
By our theorem, 
$\A$-types which are not realized by ruled surfaces must have $\A$-codimension $\ge 6$.  
This is sharp: for example, the $3$-jet $(x, y^3, x^2y)$, over which there are $\A$-orbits of codimension $6$,  
is never $\A^3$-equivalent to $3$-jets of any non-cylindrical or cylindrical ruled surfaces (Remark \ref{5jet}).  
The realizability of versal families of $\A$-types via 
families of ruled surfaces can also be verified: for each germ in Table \ref{table1}, 
an $\A_e$-versal deformation is obtained 
via deforming three invariants $\kappa_1, \tau_0, \tau_1$ appropriately (Remark \ref{versal}). 
\end{rem}

\begin{rem}\label{rem_horo}\upshape 
({\bf Conformal GA}) 
Our approach would be applicable 
to other Clifford algebras and corresponding geometries. 
For instance, Izumiya-Saji-Takahashi \cite{IST} 
classified local singularities of horospherical flat surfaces in Lorentzian space 
(conformal spherical geometry); 
a horospherical surface is described by a curve in the Lie algebra $\mathfrak{so}(3,1)$. 
Conformal Geometric Algebra may fit with this setting 
and our approach should work. 
\end{rem}

\begin{rem}\upshape 
({\bf Framed curves}) 
Take the space of dual vectors $\bD^3$ instead of $\cbU$. 
A curve $I \to \bD^3$ corresponds to a {\it framed curve}, 
which describes a $1$-parameter family of Euclidean motions of $\R^3$; 
 various geometric aspects of framed curves have 
recently been studied by e.g. Honda-Takahashi \cite{HT}. 
Since the dual Frenet formula is available for regular framed curves, 
we may rebuild the theory by using dual quaternions. 
That would be useful for singularity analysis in several topics of applied mathematics 
such as 3D-interpolation via ruled/developable surfaces, 1-parameter motions of axes in robotics,  
and so on (cf. \cite{PW, Selig}). 
\end{rem}

\subsection{Developable surfaces}  
For a non-cylindrical ruled surface, 
it is developable (the Gaussian curvature zero) if and only if 
$\kappa_1= 0$ identically (see \S 2). 
Thus two real functions $\tau_0, \tau_1$ are complete invariants of developables. 
Izumiya-Takeuchi  \cite{IT} classified generic singularities of developable surfaces rigorously,  
and Kurokawa \cite{Kurokawa} treated with a similar task for $1$-parameter families of developables. 
We generalize their results systematically using the complete invariants. 

\begin{thm}\label{thm2}
The $\A$-classification of singularities of $F$  
arising in generic at most $2$-parameter families of non-cylindrical developable surfaces 
is given as in Table \ref{table2}; In particular, 
for each $\A$-type in that table, 
the canonical expansion with the described condition is regarded as  
a normal form of the jet of developable-germ under rigid equivalence.  
\end{thm}
 
\begin{table}
{ 
$$
\begin{array}{l | l | l | l }
&\mbox{\footnotesize normal form} &\ell 
& \mbox{\footnotesize cond. at $s=s_0$}\\
\hline\hline
cE& (x, y^2, y^3) 
&1& \tau_0\not=0, \;\; \tau_1\not=0 \\
\hline 
cS_0& (x, y^2, xy^3) 
&2& \tau_1\not=0, \;\; \tau_0=0, \;\; \tau_0'\not=0  \\ 
cS_1^+&  (x, y^2, y^3(x^2+y^2)) 
&3& \tau_1\not=0, \;\; \tau_0=\tau_0'=0, \;\; \tau_0''\not=0  \\ 
cC_3^+&  (x, y^2, y^3(x^3+xy^2)) 
&4& \tau_1\not=0, \;\; \tau_0=\tau_0'=\tau_0''=0, \;\; \tau_0'''\not=0   \\
\hline 
Sw& (x, xy + 2 y^3,  xy^2 + 3 y^4) 
&2&\tau_0\not=0, \;\; \tau_1=0, \;\; \tau_1'\not=0 \\
cA_4& (x, xy + \frac{5}{2} y^4,  xy^2 + 4 y^5) 
&3&\tau_0\not=0, \;\; \tau_1=\tau_1'=0, \;\; \tau_1''\not=0\\
cA_5& (x, xy + 3 y^5,  xy^2 + 5y^6)^\dagger 
&4&\tau_0\not=0, \;\; \tau_1=\tau_1'=\tau_1''=0, \;\; \tau_1'''\not=0\\ 
 \hline 
T_1& (x, xy + y^3,  0) + o(3)
&3&\tau_0=\tau_1=0, \;\; \tau_1'\not=0\\
T_2& (x, xy,  0) + o(3) 
&4& \tau_0=\tau_1=\tau_1'=0 
\end{array}
$$
}
\caption{$\A$-types of singularities of developable surfaces. 
An exception is the type $cA_5$; 
the condition implies that the germ is {\it topologically} $\A$-equivalent to the normal form $\dagger$ 
(in this case, the striction curve $\sigma$ is topologically determinative in the sense of Ishikawa \cite{Ishikawa}).} 
\label{table2}
\end{table}

\begin{rem}\label{rem1}
\upshape 
({\bf Realization}) 
In our classification process \S 4.1, we see that 
non-cylindrical developables do not admit $\A$-types  
$$cS_1^-: (x, y^2, y^3(x^2-y^2))\;\; \mbox{and}\;\; cC_3^-: (x, y^2, y^3(x^3-xy^2))$$ 
(for the former, it was shown in \cite{Kurokawa}), 
while $cS_1^+$ and $ cC_3^+$ appear. 
Furthermore, $\tau_1\not=0$ and $\tau_0=\tau_0'=\tau_0''=0$  
if and only if the $5$-jet of $F$ is equivalent to $(x, y^2, 0)$, and 
thus, for instance,  we see that frontal singularities of  
cuspidal $S$ and $B$-types 
$$cS_*: (x, y^2, y^3(y^2+h(x,y^2))), \quad 
cB_*: (x, y^2, y^3(x^2+h(x,y^2)))$$ 
($h(x,y^2)=o(2)$) never appear in our developable surfaces. 
Similarly, since $\tau_1=0$ if and only if the $2$-jet is reduced to $(x, xy, 0)$,   
wavefronts of cuspidal beaks/lips type $A_3^\pm$ and purse/pyramid types $D_k$ never appear.  
Indeed, their $2$-jets are equivalent to $(x, 0, 0)$ and $(x^2\pm y^2, xy,0)$ respectively 
(it is obvious to see no appearance of $D_k$, for the corank of our maps $F$ is at most one). 
\end{rem}

A non-cylindrical developable surface, which is not a cone, 
is re-parametrized as the tangent developable of the {\it striction curve} $\sigma(s)$ (Lemma \ref{striction_curve2}). 
Here $\sigma(s)$ may be singular; 
recall that for a possibly singular space curve, 
its tangent developable is defined by 
the closure of the union of tangent lines at smooth points; 
indeed, it is a frontal surface, 
 see \S 2.4  (cf. Ishikawa \cite{Ishikawa99}). 
A space curve-germ is called to be of type $(m, m+\ell, m+\ell+r)$ if 
it is $\A$-equivalent to the germ  
$$x=s^m+o(m), \;\; y=s^{m+\ell}+o(m+\ell), \;\; z= s^{m+\ell+r}+o(m+\ell+r)$$ 
(the curve is called to be of {\it finite type} if $m, n, \ell <\infty$). 
A type of curve-germ is called {\it smoothly determinative} (resp. {\it topologically determinative}) 
if it determines the $\A$-type  (resp. topological $\A$-type) of the tangent developable. 
Ishikawa \cite{Ishikawa, Ishikawa99} gave the following complete characterization 
(Mond \cite{Mond89} for the case of $m=1$, i.e. smooth curves): 
\begin{itemize}
\item[(i)] smoothly determinative types are only $(1, 2, 2+r)$, $(2, 3, 4)$, $(1, 3, 4)$, $(3, 4, 5)$ and $(1, 3, 5)$; 
\item[(ii)] $(m, m+\ell, m+\ell+r)$ is topologically determinative 
if and only if $\ell, r$ are not both even, or $m=1$ and $\ell, r$ are both even. 
\end{itemize}

Using this result, 
we obtain  a complete topological $\A$-classification of singularities of 
non-cylindrical developable surfaces: 

\begin{thm}\label{thm3} {\bf (Topological classification)} 
For a non-cylindrical developable surface, 
the germ of its striction curve $\sigma(s)$ at $s=s_0$ has the type 
$$(m, m+1, m+1+r)$$
where $m-1$ and $r-1$ are orders of $\tau_1$ and $\tau_0$ at $s=s_0$, respectively: 
$$\tau_1=\tau_1'=\cdots =\tau_1^{(m-2)}=\tau_0=\tau_0'=\cdots = \tau_0^{(r-2)}=0, 
\quad \tau_1^{(m-1)}\tau_0^{(r-1)}\not=0.$$
In particular, 
topological $\A$-type of 
the germ of $F$ at singular points are completely determined by  
orders of the dual torsion $\check{\tau}=\tau_0+\varepsilon \tau_1$. 
\end{thm}

\begin{rem}\label{mond_tangent}\upshape 
Theorem \ref{thm3} is regarded as the {\it dual  version} of 
a result of Mond \cite{Mond89} and Ishikawa \cite{Ishikawa}: 
$\A$-type of the tangent developable of 
a {\it non-singular} space curve $\sigma$ with non-zero curvature is determined 
by the vanishing order of its torsion function.  
This is the case that $\sigma$ is of type $(1,2,2+r)$,  
and then  the torsion of $\sigma$ has the same order of $\tau_0$ (Lemma \ref{striction_curve2}). 
Note that in our theorem above, $\sigma(s)$ can be singular and  
non-zero curvature condition is replaced by the non-cylindrical condition. 
\end{rem}

\begin{rem}\upshape 
Table \ref{table2} is separated into three parts. 
One is the case of $\tau_1(s_0)\not=0$; 
they are the tangent developables of non-singular curves of type $(1,2,2+r)$, 
which are {\it frontal singularities} as mentioned in Remark \ref{mond_tangent}. 
The second is the case of $\tau_0(s_0)\not=0$; 
they are the tangent developables of singular curves of type 
$(2,3,4)$, $(3,4,5)$ and $(4,5,6)$, 
which are {\it wavefronts} 
-- the former two types are smoothly determinative, 
while the third one is topologically determinative, by Ishikawa's characterization. 
In the remaining part,  types $T_0$ and $T_1$ are tangent developable of curves of type 
$(2,3,4+r)\; (r \ge 1)$. 
Tangent developables of curves of types  
$(1,3,3+r), (2,4,4+r)$ etc are cylindrical at $s=s_0$. 
\end{rem}

\begin{rem}\upshape 
Not only striction curves but also several other kind of characteristic curves on a ruled surface 
can be discussed. For instance, flecnodal curves are important 
in projective differential geometry of surfaces \cite{Kabata, SKSO}. 
\end{rem}

The rest of this paper is organized as follows. 
In \S 2, we briefly review two main ingredients for non-experts in each subject -- 
the first is the algebra of dual quaternions, 
which is the most basic Geometric Algebra, and the second is about 
useful criteria for detecting $\A$-types in singularity theory of maps. 
In \S 3, we apply the $\A$-criteria to the canonical expansion of $F$ 
at singular points and prove Theorem \ref{thm1}. 
In \S 4, we proceed to the case of developable surfaces 
and prove Theorems \ref{thm2} and \ref{thm3}. 

This paper is based on the first author's  \cite{Tanaka}. 
The second author was partly supported by JSPS KAKENHI Grant Number 15K13452. 

\section{Preliminaries}
{\it Geometric Algebra} is a new look at Clifford algebras, which is nowadays recognized as a very neat tool for describing motions in Klein geometries in the context of a vast of applications to physics, mechanics and computer vision. 
In \S2.1 and \S2.2, we give a very quick summary on the geometric algebra for $3$-dimensional Euclidean motions and its application to the geometry of ruled surfaces. A good compact reference is the nineth chapter of Selig  \cite{Selig} 
(also see \cite{Guggenheimmer, IT, ICRT, PW}). 

In 2.3 and 2.4, we briefly describe some basic notions in Singularity Theory, 
which will be used in Sections 3 and 4. 
We deal with two classes of $C^\infty$ maps from a surface into $\R^3$; 
ordinary smooth maps of corank at most one, i.e. $\dim \ker df\le 1$ (Mond \cite{Mond85}) 
and frontal maps (Ishikawa \cite{Ishikawa}, Izumiya-Saji \cite{IS}). 

\subsection{Dual quaternions} 
Let $\bH$ denote the field of quaternions: 
$q=a+bi+cj+dk$. 
The conjugate of $q$ is $\bar{q}=a-bi-cj-dk$ and 
the norm is given by $|q|=\sqrt{q\bar{q}}$. 
Decompose $\bH$ into the real and the imaginary parts,  
$\bH=\R\oplus \Imag\, \bH$, 
where one identifies $bi+cj+dk \in \Imag\, \bH$ with $\bv=(b, c, d)^T \in \R^3$ 
equipped with the standard inner and exterior products. 
We write $q=a+\bv$, then the multiplication of $\bH$ is written as 
$(a+\bv)(b+\bu)=
(ab-\bv\cdot \bu) + (a\bu+b\bv +\bv \times \bu)$.  
The quaternionic unitary group 
$$\bH_1=Sp(1)=\{q \in \bH, |q|=1\}$$ 
is naturally isomorphic to $SU(2)$, that doubly covers $SO(3)$; 
indeed, $\pm q\in \bH_1$ defines the rotation $\bx \mapsto q\bx \bar{q}$. 
The Lie algebra of $\bH_1$ is just $\Imag\, \check{\bH}=\R^3$. 

Put $\bD=\R[\varepsilon]/\langle \varepsilon^2 \rangle$, 
and call it the algebra of {\it dual numbers}. 
A dual number $a+\varepsilon b$ is invertible if $a\not=0$, 
and it has a square root if $a>0$. 
The $\R$-algebra of {\it dual quaternions} is defined by  
$$\check{\bH}:= \bD^4 =\bH \otimes_\R \bD = \{\;  
\check{q}=q_0+\varepsilon q_1 \; | \;  q_0, q_1 \in \bH \;\}.$$ 
That is identified with the even Clifford algebra $C\ell^+(0,3,1)$ \cite[\S 9.3]{Selig}. 
The conjugate of $\check{q}$ is defined by $\check{q}^*:=\bar{q}_0+\varepsilon \bar{q}_1$, 
and then $\check{q}\check{q}^*=|q_0|^2+\varepsilon {\rm Re}[q_1\bar{q}_0]$. 
The Lie group of {\it unit dual quaternions} is defined by 
$$\check{\bH}_1:=\{\; \check{q} \in \check{\bH} \; | \; \check{q}\check{q}^*=1\; \}.$$
This group is isomorphic to the semi-direct product 
$\bH_1 \ltimes \Imag\, \bH=Sp(1) \ltimes \R^3$ via 
the correspondance $\check{q} \leftrightarrow (q_0, q_1\bar{q}_0)$. 
Then, $\check{\bH}_1$ doubly covers $SE(3)=SO(3) \ltimes \R^3$, 
the group of Euclidean motions of $\R^3$; the action $\check\Theta$ of $\check{\bH}_1$ on $\bx \in \R^3$ is given by 
$$1+\varepsilon \check\Theta(\check{q})\bx:=\check{q}(1+ \varepsilon \bx) \check{q}^*
=1+\varepsilon (q_0\bx \bar{q}_0+2q_1\bar{q}_0).$$ 
That is, $q_0$ and $2q_1\bar{q}_0$ express a rotation and a parallel transition, respectively. 
The Lie algebra of $\check{\bH}_1$ is canonically identified with 
the space of  {\it dual vectors} 
$$\bD^3 = \Imag\, \bH \otimes_\R \bD, 
\quad 
\check{\bv}=\bv_0+\varepsilon \bv_1 \quad (\bv_0, \bv_1 \in \Imag\, \bH=\R^3),$$
which is a $\bD$-submodule of $\check{\bH}=\bD^4$. 
The standard inner and exterior products of $\R^3$ are extended to 
$\bD$-bilinear operations on $\bD^3$;  
$$\textstyle 
\check{\bu}\cdot \check{\bv}
:= - \frac{1}{2}(\check{\bu}\check{\bv}+\check{\bv}\check{\bu}) \in \bD, 
\qquad
\textstyle
\check{\bu}\times \check{\bv}:=
\frac{1}{2}(\check{\bu}\check{\bv}-\check{\bv}\check{\bu}) \in \bD^3. 
$$
A {\it unit dual vector} means a dual vector $\check{\bv} \in \bD^3$ with 
$\check{\bv}\cdot \check{\bv}=1$, i.e. $ |\bv_0|=1$, $\bv_0\cdot \bv_1=0$ 
(indeed, it corresponds to a $2$-blade in the Clifford algebra \cite[\S 10.1]{Selig}). 
Denote the set of unit dual vectors by $\cbU$, 
which is identified with 
the space of oriented lines in $\R^3$ in the following way: 
$$\mbox{oriented lines}: \bv_0\times \bv_1+t \bv_0 \;\;  
\stackrel{1:1}{\longleftrightarrow} \;\; 
\mbox{unit dual vectors}: \bv_0+\varepsilon \bv_1.$$
This expression is very useful \cite[\S9.3]{Selig}: For instance, 
\begin{enumerate}
\item[(i)] 
a point $\ba \in \R^3$ lies on the line corresponding to 
a unit dual vector $\bv_0+\varepsilon \bv_1$ 
if and only if $\ba \times \bv_0=\bv_1$; 
\item[(ii)] 
two lines intersect perpendicularly if and only if 
corresponding unit dual vectors $\check{\bu}$ and $\check{\bv}$ satisfy that 
$\check{\bu}\cdot \check{\bv}=0$. 
\end{enumerate}

\subsection{Ruled and developable surfaces} 
Using the identification just mentioned above, 
a ruled surface is exactly described as a curve  of unit dual vectors: 
$$\cbl: I \to \cbU \subset \bD^3, \quad 
\cbl(s)=\bv_0(s)+\varepsilon \bv_1(s)$$ 
($I$ an open interval) with 
$|\bv_0(s)|=1$ and $\bv_0(s)\cdot\bv_1(s)=0$ ($s \in I$).  
Interpreting it as an object in $\R^3$,  we have a parametrization 
$$
F(s,t)=\br(s)+t\be(s) \quad (\br=\bv_0\times\bv_1,\quad\be=\bv_0). 
$$
Note that $|\be(s)|=1$ and $\br\cdot \be=0$. 
Let $R_{s}$ denote the ruling defined by $\cbl(s)$ and put 
$$R=R(\cbl):=\bigcup_{s\in I} R_s \; \subset \R^3.$$ 
Formally, $\cbl(s)$ looks like a $\bD$-version of the velocity vector of a space curve. 
That leads us to  define the {\it curvature} $\ckappa(s)$ of $\cbl$ by 
$$\ckappa(s)=\kappa_0(s)+\varepsilon\kappa_1(s)
:=\sqrt{\cbl'(s)\cdot\cbl'(s)}
=|\bv_0'|+
\varepsilon \frac{\;\; \bv_0'\cdot\bv_1'}{|\bv_0'|}
 \; \in \bD, 
$$
provided $\cbl$ is non-cylindrical, i.e., $\bv_0'(s)\not=0\; (s \in I)$.  
Here $\prime$ means $\frac{d}{ds}$. 
From now on, we assume that 
$$|\bv_0'(s)|=1$$ 
by taking $s$ to be the arc-length of $\bv_0$. 
Then, $\ckappa=1+\varepsilon \bv_0'\cdot \bv_1'$ and 
thus $\ckappa^{-1}=1-\varepsilon \bv_0'\cdot \bv_1'$.  
Put 
$$
\cbn(s) =\bn_0(s)+\varepsilon\bn_1(s):=\ckappa^{-1}\cbl'(s),
$$
and
$$
\cbt(s)={\bt}_0(s)+\varepsilon{\bt}_1(s):=\cbl(s)\times\cbn(s). 
$$ 
Then  for every $s \in I$, three dual vectors 
$\cbl(s)$, $\cbn(s)$ and $\cbt(s)$ form a basis of 
the $\bD$-module ${\rm Im}\, \check{\bH}=\bD^3$ satisfying 
$$
\begin{array}{c}
\cbl\times\cbn=\cbt,\quad
\cbt\times\cbl=\cbn,\quad
\cbn\times\cbt=\cbl,
\\
\cbl\cdot\cbn=\cbn\cdot\cbt
=\cbt\cdot\cbl
=0,\quad
\cbl\cdot\cbl=\cbn\cdot\cbn
=\cbt\cdot\cbt=1. 
\end{array}
$$
From these relations and the property (ii) of unit dual vectors mentioned before, 
we see that three lines corresponding to unit dual vectors $\cbl, \cbn, \cbt$ 
meet at one point and are mutually perpendicular; 
in particular, $\bv_0, \bn_0, \bt_0$ forms an orthonormal basis of $\R^3$. 

We define the {\it  torsion} $\ctau(s)$ of $\cbl$ by 
$$\ctau(s)=\tau_0(s)+\varepsilon\tau_1(s):=
\cbn'(s)\cdot\cbt(s) \;\; \in \bD.$$
The following theorem is classical: 

\begin{thm}\label{frenet}
{\rm (cf. Guggenheimmer \cite[\S 8.2]{Guggenheimmer}, Selig \cite[\S 9.4]{Selig})} 
Assume that $s$ is the arc-length of $\bv_0$, i.e. $\kappa_0(s)=|\bv_0'(s)|=1$. 
\begin{enumerate} 
\item {\rm (Frenet formula)}
It holds that 
\begin{equation}\label{frenet-eq}
\frac{d}{ds}
\left[\small\begin{array}{c}
\cbl(s)\\\cbn(s)\\\cbt(s)
\end{array}\right]
=
\left[\small\begin{array}{ccc}
0&\ckappa(s)&0
\\ -\ckappa(s)&0&\ctau(s)
\\0&-\ctau(s)&0
\end{array}\right]
\left[\small\begin{array}{c}
\cbl(s)\\\cbn(s)\\\cbt(s)
\end{array}\right]. 
\nonumber
\end{equation}
\item 
The dual curvature $\ckappa(s)$ and the dual torsion $\ctau(s)$ 
are complete invariants of the ruled surface $R$ up to Euclidean motions. 
That is, 
for two curves $\cbl_1$ and $\cbl_2$, 
they have the same invariants $\ckappa$ and $\ctau$ 
if and only if 
ruled surfaces $R(\cbl_1)$ and $R(\cbl_2)$ in $\R^3$ are transformed to 
each other by some Euclidean motion. 
\item 
$R(\cbl)$ is a developable surface (including a cone) if and only if $\kappa_1=0$ identically.  
In particular,  $\tau_0, \tau_1$ are complete invariants of the developable surface. 
\end{enumerate}
\end{thm}

The {\it striction curve} of a ruled surface $R$ is the curve having minimal length 
which meets all the rulings of $R$. 
Let $F(s,t)=\br(s)+t\be(s)$ be a canonical parametrization 
($\br\cdot \be=0$, $|\be|=|\be'|=1$), 
then the striction curve $\sigma(s)$ is characterized by the equation 
$\sigma'\cdot \be'=0$ (cf. \cite[p.218]{Selig}, \cite[Lemma 2.1]{IT}, \cite[\S5.3]{PW}). 
We then have the following:
 
\begin{lem} \label{striction_curve0}
For a non-cylindrical ruled surface, 
it holds that 
\begin{enumerate}
\item 
$\sigma(s)=\br(s)-(\br'(s)\cdot \be'(s))\be(s)$, 
\item  
$\sigma \times \bv_0=\bv_1$, 
$\sigma \times \bn_0=\bn_1$ and 
$\sigma \times \bt_0=\bt_1$, 
\item 
$\sigma'(s)=\tau_1(s) \bv_0(s)+\kappa_1(s)\bt_0(s)$, 
\item 
$\kappa_1 =\det(\be, \be',\br')$,   
$\tau_0=\det(\be,\be',\be'')$, 
$\tau_1 = \sigma' \cdot \be$.  
\end{enumerate}
\end{lem}

From (2) and the property (i) of unit dual vectors in \S2.1, 
it follows that $\sigma(s)$ lies on each of three lines corresponding to unit dual vectors $\cbl(s), \cbn(s), \cbt(s)$, 
that is, $\sigma(s)$ is the locus of the center of  moving orthogonal frames. 
For completeness we prove the lemma, although it is elementary. 

\

\proof 
It is easy to see (1) by differentiating $\sigma(s)=\br(s)+t(s)\be(s)$. 
We show (2). 
First, by $\cbn\cdot \cbl=0$, 
we see that 
$\bn_1\cdot \bv_0=-\bv_1\cdot \bn_0$, and similarly $ \bn_1\cdot \bt_0=-\bt_1\cdot \bn_0$. 
By the Frenet formula, 
$\bv_0'=\bn_0$, $\bt_0'=-\tau_0\bn_0$, $\bv_1'=\kappa_1\bn_0+\bn_1$ and 
$\bt_1'=-\tau_0\bn_1-\tau_1\bn_0$. 
Since $\br=\bv_0\times \bv_1$ and $\be=\bv_0$, 
it follows from (1) that  
$$\sigma=
-(\bt_1\cdot \bn_0)\bv_0-(\bv_1\cdot \bt_0)\bn_0-(\bn_1\cdot \bv_0)\bt_0.$$ 
Thus 
$\sigma \times \bv_0=
-(\bv_1\cdot \bt_0)\bn_0\times \bv_0-(\bn_1\cdot \bv_0)\bt_0\times \bv_0
=(\bv_1\cdot \bt_0)\bt_0+(\bv_1\cdot \bn_0)\bn_0=\bv_1$, for $\bv_1\cdot \bv_0=0$. 
That yields (2).  
Differentiating the first one of (2),  
$$0=(\sigma\times \bv_0)'-\bv_1'=
(\sigma'\times \bv_0+\sigma \times \bn_0)-(\kappa_1\bn_0+\bn_1)
=\sigma'\times \bv_0-\kappa_1\bn_0$$ 
and similarly $\sigma'\times \bt_0+\tau_1\bn_0=0$. 
Substitute $\sigma'=a\bv_0+b\bn_0+c\bt_0$ for those equalities, 
we obtain $a=\tau_1$, $b=0$, $c=\kappa_1$, that is (3). 
Finally, (4) is easy, e.g.,  $\kappa_1 = \bv_0'\cdot \bv_1'=\be'\cdot (\br'\times \be)
=\det(\be, \be',\br')$. 
\qed

\begin{lem}\label{striction_curve}
 {\rm (Izumiya et al \cite[Lemma 2.2]{IT}, \cite[\S 1]{ICRT})}
For a non-cylindrical ruled surface, 
$F$ is singular at $(s_0, t_0)$ if and only if $\kappa_1(s_0)=0$ and 
$t_0=-\br'(s_0)\cdot \be'(s_0)$. 
The singular value $F(s_0, t_0)$ is the point $\sigma(s_0)$ 
where the curve $\sigma(s)$ is tangent to the ruling $R_{s_0}$ or 
$\sigma'(s_0)=0$. 
\end{lem}

\proof 
$\frac{\rd F}{\rd s}(s_0) \times \frac{\rd F}{\rd t}(s_0)
=(\br'(s_0)+t_0\be'(s_0))\times \be(s_0)=0$ 
$\Leftrightarrow$ $\br'(s_0)=\alpha \be(s_0)-t_0\be'(s_0)$ 
for some $\alpha\not=0$ 
$\Leftrightarrow$ 
$\det(\be(s_0),\be'(s_0),\br'(s_0))=0$ and $t_0=-\br'(s_0)\cdot \be'(s_0)$. 
The second claim follows from (3) in Lemma \ref{striction_curve0}. 
\qed

\

In case of  $\kappa_1=0$ identically, Lemmas \ref{striction_curve0} and \ref{striction_curve} imply that 
singular points of $F$ form a non-singular curve 
$s \mapsto (s, -\br'(s)\cdot \be'(s)) \in I \times \R$ 
and the image of this curve is just the striction curve $\sigma(s)$. 
Note that 
$\sigma(s)$ is a non-singular space curve, if $\tau_1\not=0$; 
especially, 
$F$ is written by $\sigma(s)+\tilde{t}\sigma'(s)$ with 
$\tilde{t}=(t+\br'(s)\cdot \be'(s))/\tau_1$.

\begin{lem} {\rm (Izumiya et al \cite[\S 1]{ICRT})} \label{striction_curve2}
A non-cylindrical developable surface, which is not a cone, 
is re-parametrized as the tangent developable of the striction curve $\sigma(s)$. 
The curve $\sigma$ is non-singular whenever $\tau_1\not=0$, and then 
the curvature $\kappa_\sigma$ and the torsion $\tau_\sigma$ of $\sigma$ are given respectively by  
$$\kappa_\sigma=\frac{|\sigma'\times \sigma''|}{|\sigma'|^3}=\frac{1}{\tau_1}, \qquad 
\tau_\sigma=
\frac{\det(\sigma', \sigma'', \sigma''')}{|\sigma'\times \sigma''|^2}=\frac{\tau_0}{\tau_1}.$$
\end{lem}

\subsection{$\A$-classification of map-germs} 
A {\it singular} point of $f: M \to N$ between manifolds means a point $p \in M$ where 
$df_p$ is neither injective nor surjective 
(then $f(p) \in N$ is called a {\it singular value} of $f$);  
we denote by $S(f)\subset M$ the set of singular points of $f$. 
Two maps $\tilde{f}: U \to N$ and $\tilde{g}: V \to N$ 
on neighborhoods $U$ and $V$ of  $p\in M$ 
define the same {\it map-germ at $p$} 
if there is a neighborbood $W\subset U\cap V$ of $p$ 
so that $\tilde{f}|_W\equiv \tilde{g}|_W$;  
a {\it map-germ at $p$} is an equivalence class of maps under this relation, 
denoted by $f: (M, p) \to (N, f(p))$. 
Two map-germs at $p$ have the same {\it $k$-jet} 
if they have the same Taylor polynomials at $p$ of order $k$ in some local coordinates; 
a $k$-jet is such an equivalence class of map-germs, denoted by $j^kf(p)$. 
Two germs $f: (M, p) \to (N, q)$ and $g: (M', p') \to (N', q')$ are {\it $\A$-equivalent} if 
they commute each other 
via diffeomorphism-germs $\sigma$ and $\tau$:  
$$
\xymatrix{
(M, p) \ar[r]^f  \ar[d]_{\sigma}^{\simeq} &  (N, q) \ar[d]^{\tau}_{\simeq}\\
(M', p')\ar[r]_g &  (N', q')
}
$$
For simplicity, 
we consider map-germs $(\R^m,0) \to (\R^n,0)$ and the $\A$-equivalence by 
the action of diffeomorphisms $\sigma$ and $\tau$ preserving the origins. 
At the $k$-jet level, {\it $\A^k$-equivalence} is defined. 
A germ $f:(\R^m,0) \to (\R^n,0)$ is said to be {\it $k$-$\A$-determined} if 
any germs $g:(\R^m,0) \to (\R^n,0)$ with $j^kg(0)=j^kf(0)$ is $\A$-equivalent to $f$; 
such germs are collectively referred to as  {\it finitely $\A$-determined} germs. 
For instance,  the germ $(x, y^2, xy)$ is $2$-determined. 
Let $J^k(m,n)$ be the {\it jet space} consisting of all $k$ jets of $(\R^m,0) \to (\R^n,0)$, 
which is identified with 
the affine space of Taylor coefficients of order $r$ ($1\le r\le k$) 
in a fixed system of local coordinates. 
The codimension of the $\A$-orbit of a germ $f$ in 
the space of all map-germs $(\R^m,0) \to (\R^n,0)$ 
is called the {\it $\A$-codimension} of $f$; 
the $\A$-codimension of $f$ is finite if and only if $f$ is finitely $\A$-determined (see e.g. \cite{GG}). 

Thanks to finite determinacy, 
the process of $\A$-classification is reduced to a finite dimensional problem: 
we stratify $J^k(m,n)$ invariantly under the $\A^k$-equivalence step by step 
from low order $k$ and low codimension. 
For instance, 
using several determinacy criteria, 
$\A$-classification of map-germs $(\R^2, 0) \to (\R^3, 0)$ up to certain codimension 
has been established in Mond \cite{Mond82, Mond85}. 
In \S 3, we will follow Mond's classification process. 

Furthermore, in Mond \cite{Mond82, Mond89},  
a special class of map-germs  $(\R^2, 0) \to (\R^3, 0)$ is considered. 
A map germ $f: (\R^2,0) \to (\R^3, 0)$ is of class {\it CE} (i.e. cuspidal edge), 
if $\rank df(0)=1$ and the singular point set $S(f)$ is non-singular. 
A germ $f$ in CE is {\it  $k$-$\A$-determined in CE} 
if any germ $g$ in CE with the same $k$-jet 
as $j^kf(0)$ is $\A$-equivalent to $f$. 
In \S 4, we will use the following criteria of determinacy in CE  \cite[Lem.1.1, Prop.1.2]{Mond89}. 

\begin{prop}[Mond \cite{Mond89}] \label{CE} It holds that \\
$i)$ If $f \in CE$ and $j^2f(0)=(x,y^2,0)$, then $f$ is $\A$-equivalent to the germ $g(x,y)=(x,y^2,y^3p(x,y^2))$ for some smooth function $p(u,v)$; \\ 
$ ii)$ $f(x,y)=(x,y^2, y^3)$ is $3$-determined in CE; \\ 
$ iii)$  $f(x,y)=(x, y^2, yp(x, y^2))$ and $g(x,y)=(x, y^2, yq(x, y^2))$ are $\A$-equivalent if and only if $\tilde{f}(x,y)=(x, y^2, y^3p(x, y^2))$ and $\tilde{g}(x,y)=(x, y^2, y^3q(x, y^2))$ are $\A$-equivalent. In particular, $f$ is $(k-2)$-determined if and only if $\tilde{f}$ is $k$-determined in CE. \end{prop}

\subsection{Singularities of frontal surfaces} 
There is a special class of surfaces, called {\it frontal surfaces}. 
Let $ST^*\R^3$ be the spherical cotangent bundle with respect to the standard metric of $\R^3$ 
equipped with the standard contact structure. Let $U$ be an open set of $\R^2$. 
A map  $\iota: U \to ST^*\R^3$ is called {\it isotropic} if  it satisfies that 
the image $d\iota(T_pU)$ is contained in the contact plane $K_{\iota(p)}$ for any $p \in U$. 
A {\it frontal map} 
is the composed map $f=\pi\circ \iota: U \to \R^3$ of an isotropic map $\iota$ and the projection 
$\pi: ST^*\R^3\to \R^3$. 
The image (possibly singular) surface is called to be {\it frontal}. 
An isotropic immersion $\iota$ is usually called a {\it Lagrange immersion}, 
and $\pi\circ \iota$ and its image are called a {\it Lagrange map} and a {\it wavefront}, respectively. 
Let $f: U \to \R^3$ be a frontal map with $\nu: U \to S^2$ so that 
$\iota=(f, \nu): U \to ST^*\R^3=\R^3\times S^2$ is an isotropic map.  
We identifies $T\R^3\simeq T^*\R^3$ using the standard metric, 
then the unit vector $\nu$ is always orthogonal to the subspace $df(T_pU)$ at any $p \in U$. 
Let $x, y$ be coordinates of $U$ and 
put $\lambda(x,y)=\det \left[\frac{\rd f}{\rd x}, \frac{\rd f}{\rd y}, \nu \right](x,y)$; 
then  the singular point set $S(f)$ is defined by $\lambda(x,y)=0$. 
If $d\lambda(p)\not=0$, then $p$ is called a {\it non-degenerate} singular point. 
In particular, if $p$ is non-degenerate and $\mbox{rank}\, df_p=1$, 
the germ $f$ at $p$ is of class CE. 

For a developable surface with $\be\times \be'\not=0$,  
set $f: U \to \R^3$ to be $f(s,t):=\br(s)+t\be(s)$. 
Then $f$ is a frontal map; 
in fact, it suffices to put $\nu=\be\times \be'/|\be\times \be'|$ 
(then $\frac{\rd f}{\rd t}\cdot \nu=\be\cdot \nu=0$ and 
$\frac{\rd f}{\rd s} \cdot \nu=(\br'+t\be')\cdot \nu=\det(\br', \be, \be')=0$). 
Note that any singularities of $f$ are non-degenerate and have corank one 
(see the comment before Lemma \ref{striction_curve2}). 
There are two cases: 

If  $\iota=(f, \nu)$ is singular, 
then it is easy to see that the $2$-jet of $f$ is $\A^2$-equivalent to $(x,y^2,0)$, and hence 
Mond's criteria for map-germs of class CE (Proposition \ref{CE}) can be applied. 

If $\iota$ is non-singular, i.e. $\iota$ is a Legendre immersion, 
then the $2$-jet is equivalent to $(x, xy,0)$, and thus Proposition \ref{CE} is useless. 
In this case, we employ the Legendre singularity theory. 
There are known useful criteria of \cite{IS} 
(precisely saying, the topological type $cA_5$ is not dealt in \cite{IS} but 
the same argument as in Appendix of \cite{IS} works as well): 

\begin{prop} {\rm (Izumiya-Saji  \cite[Theorem 8.1]{IS}) }
Let $f: U \to \R^3$ be a Legendre map, and $p$ 
a non-degenerate singular point with $\mbox{rank}\, df_p=1$. 
Let $\eta$ be an arbitrary vector field around $p$ 
so that $\eta(q)$ spans $\ker df_q$ at any $q \in S(f)$. 
Then $f$ is $\A$-equivalent to $cE$,  $Sw$,  $cA_4$ or $cA_5$ 
if the following condition holds: 
$$
\begin{array}{l | l}
cE & \eta\lambda(p)\not=0\\
Sw &  \eta\lambda(p)=0, \; \eta\eta\lambda(p)\not=0 \\
cA_4 &  \eta\lambda(p)=\eta\eta\lambda(p)=0, \; \eta\eta\eta\lambda(p)\not=0 \\
cA_5 &  \eta\lambda(p)=\eta\eta\lambda(p)=\eta\eta\eta\lambda(p)=0, \; 
\eta\eta\eta\eta\lambda(p)\not=0. 
\end{array}
$$
\end{prop}

Through the theory of frontal maps and generating functions, 
Ishikawa \cite{Ishikawa, Ishikawa99} 
showed that the tangent developable of a curve of type 
$(m, m+\ell, m+\ell+r)$ has a parametrization $F: (\R^2, 0) \to (\R^3,0)$ defined by 
\begin{eqnarray*}
x&=&t\\
y&=&s^{m+\ell}+s^{m+\ell+1}\varphi(s)+t(s^\ell+s^{\ell+1}\phi(s)),\\
z&=&(\ell+r)(m+\ell+r)\int_0^s u^r\frac{\rd y(u,t)}{\rd u} du\\
&=& (\ell+r)(m+\ell)s^{m+\ell+r}+\cdots + t(\ell(m+\ell+r)s^{\ell+r}+\cdots)
\end{eqnarray*}
with some $C^\infty$ functions $\varphi(s)$ and $\phi(s)$. 
These two function must be related to invariants $\tau_0$ and $\tau_1$. 
It is also shown \cite[Thm 2.1]{Ishikawa} 
that the topological type of the tangent developable of a space curve is determined  
by type $(m, m+\ell, m+\ell+r)$ of the curve, unless both $\ell, r$ are even, 
as mentioned in Introduction.  

\section{Singularities of ruled surfaces}
In this section, we prove Theorem \ref{thm1} (2); 
first we give a certain stratification of the jet space of triples of functions 
$(\kappa_1, \tau_0, \tau_1)$, and 
then discuss a variant of Thom's transversality theorem. 

\subsection{Dual Bouquet formula} 
Consider a curve $\cbl: I \to \bD^3$, $\cbl(s)=\bv_0(s)+\varepsilon \bv_1(s)$,  
with $\cbl\cdot \cbl=1$ and $|\bv_0'(s)|=1$ as in \S 2.2. 
We are concerned with the germ of $\cbl$ at the origin ($s_0=0$).  
Throughout this section, 
let $\ckappa, \ctau, \ckappa',  \ctau',  \cdots$ denote 
their values at $s=0$ for short, e.g. $\ckappa'=\ckappa'(0)$, 
unless specifically mentioned. 

By iterated uses of the Frenet formula (Theorem \ref{frenet} (1)), 
we obtain 
the ``Bouquet formula" of the curve in $\bD^3$ at $s=0$; 
$$
\cbl(s)=\sum_{n=0}^r
\frac{\cbl^{(n)}(0)}{n!}s^n + o(r) \;\; \in \bD^3
$$
with 
{\small 
\begin{eqnarray*}
\cbl'(0)&=&\ckappa\,\cbn(0),\\ 
\cbl''(0)&=&-\ckappa^2\,\cbl(0)+\ckappa'\,\cbn(0)+\ckappa\ctau\,\cbt(0),\\ 
\cbl^{(3)}(0)&=&-3\ckappa\ckappa'\,\cbl(0)
+(\ckappa''-\ckappa^3-\ckappa\ctau^2)\,\cbn(0)
+(2\ckappa'\ctau+\ckappa\ctau')\,\cbt(0), \\
\cbl^{(4)}(0)&=&(\ckappa^4+\ckappa^2\ctau^2-4\ckappa\ckappa'')\cbl(0)+(\ckappa^{(3)}-6\ckappa^2\ckappa'-3\ckappa'\ctau^2-3\ckappa\ctau\ctau')\cbn(0)
\\&&+(3\ckappa''\ctau+3\ckappa'\ctau'-\ckappa^3\ctau+\ckappa\ctau''-\ckappa\ctau^3)\cbt(0), \\
\cbl^{(5)}(0)
&=&(10\ckappa^3\ckappa'+5\ckappa\ckappa'\ctau^2+5\ckappa^2\ctau\ctau'-5\ckappa\ckappa^{(3)})\cbl(0)
+(\ckappa^{(4)}-6\ckappa^2\ckappa''-6\ckappa''\ctau^2\\
&&-12\ckappa'\ctau\ctau'-3\ckappa(\ctau')^2-4\ckappa\ctau\ctau''+\ckappa^3\ctau^2+\ckappa\ctau^4)\cbn(0)
+(4\ckappa^{(3)}\ctau+6\ckappa''\ctau'\\
&&+3\ckappa'\ctau''-9\ckappa^2\ckappa'\ctau-\ckappa^3\ctau'+\ckappa'\ctau''+\ckappa\ctau^{(3)}-4\ckappa'\ctau^3-6\ckappa\ctau^2\ctau')\cbt(0),
\end{eqnarray*}
}
and so on. 
A similar (but more naive) expansion written by Pl\"ucker coordinates, instead of dual quaternions, 
can be found in a classical literature, Hlavat\'y \cite{Hlavaty}. 

Since  dual vectors $\{\cbl(0), \cbn(0), \cbt(0)\}$ form a $\bD$-basis of $\Imag\, \check{\bH}=\bD^3$,  
we write 
$$\textstyle 
\cbl(s)=[\, \cbl(0),\cbn(0),\cbt(0)\,]\; \check{\bw}(s)$$
where 
$$\textstyle 
 \check{\bw}(s)
=[1-\frac{1}{2}\ckappa^2 s +\cdots, \ckappa s + \frac{1}{2}\ckappa's + \cdots, \frac{1}{2}\ckappa\ctau s^2+\cdots]^T \in \bD^3.$$  
Recall that 
three oriented lines in $\R^3$ determined by unit dual vectors 
$\cbl(0), \cbn(0), \cbt(0)$ meet at one point, 
which is nothing but the striction point $\sigma(0)$, 
as mentioned just after Lemma \ref{striction_curve0}. 
By an Euclidean motion, 
the three lines can be transformed to standard coordinate axises of $\R^3$, i.e., 
$\bv_0(0), \bn_0(0), \bt_0(0)$ form 
the standard basis $i, j, k$ of $\Imag\, \bH = \R^3$, respectively, 
and $\bv_1(0)=\bn_1(0)=\bt_1(0)=0 \in \R^3$. 
Namely, we may assume that 
the $3\times 3$ matrix (with entries in $\bD$) $[\,\cbl(0),\cbn(0),\cbt(0)\,]$ is the identity,  
so $\cbl(s)=\check{\bw}(s)$. 
Then, 
$$\cbl(s)=\bv_0(s)+\varepsilon \bv_1(s)=
\left[\begin{array}{c}1\\s\\0
\end{array}\right]
+\varepsilon  
\left[\begin{array}{c}
0\\  \kappa_1 s\\0
\end{array}\right]+o(1). 
$$
At a point $(0, t_0) \in I \times \R$,  the Taylor expansion  of 
$F(s,t)=\bv_0(s)\times \bv_1(s)+t\bv_0(s)$ 
is immediately obtained; 
in particular, $F(0,t_0)=[t_0,0,0]^T$ and 
$$
dF(0, t_0)= 
\left[\begin{array}{cc}
0&1\\t_0&0\\ \kappa_1&0
\end{array}\right].  
$$
This gives an alternative proof of Lemma \ref{striction_curve}: 
$F$ is singular at $(0, t_0)$ if and only if $\kappa_1(0)=t_0=0$
($t_0=0$ means that the point is just the striction point $\sigma(0)$ lying on the ruling). 
Assume that $F$ is singular at the origin. Then 
we obtain a {\it canonical Taylor expansion} of $F$: 
\begin{eqnarray}\label{eq1}
&&F(s,t)=
\textstyle
\left(t-\frac{1}{2}ts^2+\frac{\tau_1}{2} s^3,\ 
ts-\frac{\tau_1}{2}s^2-\frac{2\tau_0\kappa_1'+\tau_1'}{6}s^3,\ 
\frac{\kappa_1'}{2}s^2+\frac{\tau_0}{2}t s^2+\frac{\kappa_1''-2\tau_0\tau_1}{6}s^3\right)+o(3).
\end{eqnarray}

\begin{rem}\label{truncation}\upshape 
{\bf (Truncated polynomial maps)} 
Let $F(s,t)$ be as in (\ref{eq1}), 
and set $\bar{F}(s, t)=(\bar{\bv}_0(s) \times \bar{\bv}_1(s)) + t \bar{\bv}_0(s)$ 
to be a polynomial map of order $k$ with $j^k\bar{F}(0)=j^kF(0)$. 
Denote by $\bar{s}$ the arc-length of the curve $\bar{\bv}_0(s)$, then $\bar{s}:=s+o(k)$, 
and thus $k$-jets at $0$ of the dual curvature and the dual torsion do not change 
from those of $F$. 
That gives examples of polynomial ruled surfaces  
with prescribed $k$-jets of $\check{\kappa}$ and $\check{\tau}$ at a point. 
\end{rem}

\subsection{Recognition of singularity types} 
Our task is to find appropriate diffeomorphism-germs of the source and the target 
for reducing jets of $F(s,t)$ to normal forms in $\A$-classification 
step by step; for such computations, we have used the software {\it Mathematca}. 

Let $(X, Y, Z)$ be the coordinates of the target $\R^3$. 
Below,  $\kappa_1, \kappa_1', \cdots$ denote   
their values at $s=0$ unless specifically mentioned. 
From now on, 
assume that $\kappa_1(=\kappa_1(0))=0$. 
Put $y=s$ and $x=t-\frac{1}{2}ts^2+\frac{\tau_1}{2} s^3+\cdots$, 
the first component of $F$ in the form (\ref{eq1}) above. 
With this new coordinates $(x, y)$ of the source $\R^2$, we set
\begin{eqnarray}\label{2}
f(x,y)&:=&F(y, t(x, y))=(x, f_2(x,y), f_3(x,y)) \\
&=&\textstyle
\left( x,\ x y-\frac{1 }{2}\tau_{1} y^2-\frac{1 }{6}\tau_{1}' y^3,
\ \frac{1}{2}\kappa_1' y^2+\frac{1}{2}\tau_{0} x y^2 
+ \frac{1}{6 }(\kappa_{1}''-2 \tau_{0} \tau_{1})y^3\right)
+o(3).\nonumber
\end{eqnarray}
Note that $f(x,y)$ is still of the form $\tilde{\br}(y)+x\tilde{\be}(y)$. 
Now, we apply to this germ $f(x,y)$ the recognition trees in Mond's classification \cite[Figs.1, 2]{Mond85}.   
Below, $S_k^\pm$, $B_k^\pm$, $C_k^\pm$, $H_k$ and $F_4$ denote  Mond's notations of 
$\A$-simple germs  \cite{Mond85}.

\

\t
$\bullet\;\; ${\bf $2$-jet}: 
Crosscap $S_0$ is $2$-determined, thus it follows from (\ref{2}) that 
$$f\sim_\A S_0: (x,xy, y^2) \;\;\; \Longleftrightarrow \;\;\; 
\kappa_1=0, \quad 
\kappa_1' \neq 0.$$
Let $\kappa_1'=0$. Then 
the $2$-jet is equivalent to either of  $(x, xy,0)$ or $(x, y^2,0)$, 
according to whether $\tau_1= 0$ or not. 
We compute the second and third component of $f$ as 
\begin{eqnarray*}
f_2&=& \textstyle 
x y -\frac{1}{2} \tau_1 y^2 - \frac{1}{6}\tau_1' y^3
 + \frac{1}{24} ((8 - 4\tau_0^2) x y^3 +( - 5\tau_1  + 
3 \tau_0^2 \tau_1  - 3 \tau_0 \kappa_1''   -  \tau_1'' )y^4 ) 
\\
&&\textstyle
+ \frac{1}{120} (-15 \tau_0 \tau_0' x y^4  + 
(12  \tau_0 \tau_0'\tau_1  - 9  \tau_1' + 
6  \tau_0^2 \tau_1' - 6  \tau_0' \kappa_1'' - 
4  \tau_0 \kappa_1^{(3)} - \tau_1^{(3)}) y^5)
\\
&&\textstyle
+o(5), \\
f_3&=&
 \textstyle
\frac{1}{6} (3 \tau_0 x y^2  +( \kappa_1''- 2 \tau_0 \tau_1)y^3) + 
\frac{1}{24} (4  \tau_0' x y^3+( - 3\tau_1 \tau_0' - 
     3 \tau_0 \tau_1'+  \kappa_1^{(3)})y^4)\\
&&\textstyle
 + \frac{1}{120} ((25  \tau_0 - 5 \tau_0^3+5 \tau_0'')xy^4 
+(- 16 \tau_0 \tau_1 + 4  \tau_0^3 \tau_1
 - 6  \tau_0' tau_1'- 
6  \tau_0^2 \kappa_1''  \\
&&\textstyle
- 4 y^5 \tau_1 \tau_0''
- 4 y^5 \tau_0 \tau_1'' +\kappa_1^{(4)})y^5)+o(5). 
\end{eqnarray*}

\

\t
$\bullet\;\; ${\bf $3$-jet}: 
Let $\kappa_1=\kappa_1'=0$ and  $\tau_1\not= 0$.  
First, let us remove the term $xy$ from $f_2$; 
take $\bar{x}=x$ and $\bar{y}=y-\frac{1}{\tau_1}x$, 
then we see that 
\begin{equation}\textstyle 
j^3f(0) \sim 
(x, y^2+\frac{\tau_1'}{\tau_1^3} x^2y+\frac{\tau_1'}{3\tau_1} y^3, 
\kappa_1''x^2y+\frac{1}{3}\tau_1^2(\kappa_1''-2\tau_0\tau_1)y^3). 
\end{equation}
The first two components can be transformed to $(x,y^2)$ 
by a coordinate change of $(x,y)$ with identical linear part and by a target coordinate change of $(X,Y)$, 
since the plane-to-plane germ $(x,y^2)$ is  $2$-determined (stable germ). 
Hence $j^3f(0)$ is equivalent to one of the following: 
\begin{equation}\label{3-jet1}
\left\{
\begin{array}{lll}
(x, y^2, y^3\pm x^2 y)
& \kappa_1''( \kappa_1''-2 \tau_0\tau _1)\gtrless 0, \ \tau_1\neq 0
&\cdots S_1^\pm
\\(x, y^2, y^3)
& \kappa_1'' = 0,\ \tau_0\tau_1\neq 0
&\cdots S
\\(x, y^2, x^2 y)
& \kappa_1''\ =2 \tau_0\tau_1\neq 0,
&\cdots B
\\(x, y^2, 0)
& \kappa_1''=\tau_0=0, \ \tau_1\neq 0
&\cdots C
\end{array}
\right. 
\end{equation}
Note that $S_1^\pm$ is $3$-determined, thus this case is clarified. 

Let $\tau_1= 0$. Then from (\ref{2}), we have 
$$\textstyle
j^3f(0)\sim
\left( x,\ x y-\frac{1 }{6}\tau_{1}' y^3,\frac{1}{2}\tau_{0} x y^2 + \frac{1}{6 }\kappa_{1}'' y^3\right). 
$$
In the same way as above, $j^3f(0)$ is reduced to one of the following: 
\begin{equation}\label{3-jet2}
\left\{
\begin{array}{lll}
(x, xy, y^3)& 
\kappa_1''\neq  0,\ \tau_1= 0
&\cdots H
\\(x, xy+y^3, x y^2)& 
\kappa_1'' =\tau_1= 0,\ \tau_0\tau_1'\neq 0
&\cdots P
\\(x, xy, xy^2)& 
\kappa_1''=\tau_1= \tau_1'=0, \ \tau_0\neq 0
&
\\(x, xy+y^3, 0)& 
\kappa_1''=\tau_0=\tau_1= 0,\ \tau_1'\neq 0
&
\\(x, xy, 0)& 
\kappa_1''=\tau_0=\tau_1= \tau_1' =0
&
\end{array}
\right.
\end{equation}
Each of last three types has codimension $\ge 6$, so we omit them here. 
Below, 
for types $S, B, \cdots, P$ in (\ref{3-jet1}) and (\ref{3-jet2}), 
we detect $\A$-types with codimension $\le 5$ by checking higher jets and the determinacy. 

\

\t
$\bullet \;\; ${\bf $S$-type}: 
Let $\kappa_1=\kappa_1' = 0$ and $\tau_0\tau_1\neq 0$. 
Then a computation shows that 

$$
\begin{array}{rcl}
\textstyle
\kappa_1''=0
&\textstyle
\quad\Longrightarrow\quad
&
j^4 f(0)
\sim 
\left(x, y^2, y^3-\frac{\kappa_1^{(3)}}{2\tau_0\tau_1^4}x^3 y \right)
\\
\textstyle
\kappa_1''=\kappa_1^{(3)}=0
&\textstyle
\quad\Longrightarrow\quad
&
j^5 f(0)
\sim 
\left(x, y^2, y^3-\frac{\kappa_1^{(4)}}{8\tau_0\tau_1^5}x^4 y
\right)
\\
\textstyle
\kappa_1''=\kappa_1^{(3)}=\kappa_1^{(4)}=0
&\textstyle
\quad\Longrightarrow\quad
&
j^6 f(0)
\sim 
\left(x, y^2, y^3-\frac{\kappa_1^{(5)}}{40\tau_0\tau_1^6}x^5 y
\right). 
\end{array}
$$
Note that $S_k$ is $(k+2)$-determined (its codimension is $k+2$), thus 
$f$ is of type $S_k^\pm$ $(k=2,3,4)$
if and only if 
$\kappa_1=\kappa_1'=\cdots=\kappa_1^{(k)}=0$
and 
$\kappa_1^{(k+1)}\tau_0\tau_1 \lessgtr 0$ (seemingly, it is so for any $k$). 

\

\t
$\bullet \;\; ${\bf $B$-type}: 
Let  $\kappa_1=\kappa_1'=\kappa_1''-2 \tau_0\tau_1=0$ and $ \kappa_1''\neq 0$. 
Then it would be $\A$-equivalent to $B_k$-type \cite[4.1:17, Table 3]{Mond85}. For instance, 
$$j^5f(0)\sim (x, y^2, x^2 y+b_2 y^5)$$ 
with  
$$
\begin{array}{rl}
b_2=&48\tau_0^2 \tau_1^2(\tau_0^2-2)-20 (\tau_0^2 (\tau_1')^2+\tau_1^2 (\tau_0')^2)-56 \tau_0 \tau_1 \tau_0' \tau_1'-24\tau_0 \tau_1(\tau_0\tau_1''+\tau_1 \tau_0'')\\
&+20 \kappa_1^{(3)}(\tau_0 \tau_1'+\tau_1 \tau_0')-5 (\kappa_1^{(3)})^2+6 \kappa_1^{(4)} \tau_0 \tau_1. 
\end{array}
$$
Since $B_2$ is $5$-determined, 
$$f\sim_\A B_2^\pm:(x,y^2,x^2y\pm y^5) 
\;\;  \Longleftrightarrow \;\; b_2\gtrless 0.$$ 
Let $b_3$ be the coefficient of $y^7$ in the last component of $j^7f(0)$, 
which is written as a polynomial in derivatives of invariants at $s=0$, then 
$B_3^\pm: (x,y^2,x^2y\pm y^7)$ is detected by the condition that $b_2=0$ and $b_3\not=0$. 
Here $B_3$ is of codimension $5$.

\

\t
$\bullet \;\; ${\bf $C$-type}: 
Let $\kappa_1=\kappa_1'=\kappa_1''=\tau_0=0,\ \tau_1\neq 0$. 
Through 
$$\textstyle \psi(X,Y,Z)=\left(\frac{1}{\tau_1}X,\ Y,\ 
\frac{1}{\tau_1}\left(Z-a Y^2
-b X^2 Y\right)\right)$$ 
with $a=\frac{1}{4}(\kappa_1^{(3)}-3\tau_1\tau_0')$,  
$b=\frac{3}{2\tau_1^2}(\kappa_1^{(3)}-\tau_1\tau_0')$,  
we see that 
$$\textstyle
j^4 f(0)\sim 
\left(x,\ y^2,\ 
\kappa_1^{(3)}x^3 y
+(\kappa_1^{(3)}-2 \tau_1\tau_0')x y^3\right).$$ 
Since $C_3$ is $4$-determined (of codimension $5$), 
$$f \sim_\A C_3^\pm: (x, y^2, x y^3\pm x^3 y)
\;\;  \Longleftrightarrow \;\; 
\kappa_1^{(3)}(\kappa_1^{(3)}-2\tau_1\tau_0')\gtrless 0.$$ 

\

\t
$\bullet \;\; ${\bf $H$-type}: 
Let $\kappa_1=\kappa_1'=\tau_1= 0$ and $\kappa_1''\neq  0$. 
Then it would be $\A$-equivalent to $H_k$-type \cite[4.2.1:2]{Mond85}. 
A lengthy computation shows that 
$$
j^5 f(0)\sim 
\left(x,\ 
x y + h_2 y^5,\ 
y^3
\right)
$$
with 
$$
\begin{array}{rl}
h_2=
&\hspace{-5pt}
-15  \tau_0^2 (\tau_1')^3
-24  \tau_0' (\tau_1')^2 \kappa_1''
-36  \tau_1' (\kappa_1'')^2
-15  \tau_0^2 \tau_1' (\kappa_1'')^2
-24  \tau_0' (\kappa_1'')^3
-21  \tau_0 \tau_1' \kappa_1'' \tau_1''
\\&+20  \tau_0 (\tau_1')^2 \kappa_1^{(3)}
- \tau_0 (\kappa_1'')^2 \kappa_1^{(3)}
+5  \kappa_1'' \tau_1'' \kappa_1^{(3)}
-5  \tau_1' (\kappa_1^{(3)})^2
-4  (\kappa_1'')^2 \tau_1^{(3)}
+4  \tau_1' \kappa_1'' \kappa_1^{(4)}. 
\end{array}
$$
Since $H_2$ is $5$-determined, 
$$f\sim_\A H_2^\pm:(x,xy\pm y^5, y^3) 
\;\;  \Longleftrightarrow \;\; h_2\gtrless 0.$$ 
Let $h_3$ be the coefficient of $y^8$ in the middle component of $j^8f(0)$, 
then $H_3:(x,xy+ y^8, y^3)$ is detected by $h_2=0$ and $h_3\not=0$ 
($H_3$ is of codimension $5$). 

\

\t
$\bullet \;\; ${\bf $P$-type}: 
Let $\kappa_1=\kappa_1'=\kappa_1'' =\tau_1= 0$ and $\tau_0\tau_1'\neq 0$. 
Then we see that there is a polynomial $p_4$ in derivatives of $\kappa_1, \tau_0, \tau_1$ so that 
$$f\sim_\A 
P_3: (x,\ x y + y^3,\ x y^2 +p_4 y^4)$$ 
for  $p_4\not=0, \frac{1}{2}, 1, \frac{3}{2}$   \cite[\S 4.2]{Mond85}.

\begin{rem}\label{C_kF_4}\upshape 
{\bf (Characterization of $C_k$ and $F_4$)} 
Among $\A$-simple germs in Mond \cite{Mond85}, 
$S_k^\pm$, $B_k^\pm$ and $H_k$ have been discussed above, 
so there remain $C_k\; (k\ge 4)$ and $F_4$. 
Consider $C$-type above and think of the next to $C_3$-type. Namely, 
suppose that $\kappa_1^{(3)}(\kappa_1^{(3)}-2\tau_1\tau_0')=0$. 

$\bullet \;\; $
If $\kappa_1^{(3)}=0$ and $\tau_1\tau_0'\neq 0$, then $j^4 f(0)\sim (x,y^2,xy^3)$ and 
$$
\begin{array}{rcl}
\textstyle
\kappa_1^{(3)}=0
&\textstyle
\quad\Longrightarrow\quad
&
j^5 f(0)
\sim 
\left(x, y^2, x y^3-\frac{\kappa_1^{(4)}}{8\tau_0'\tau_1^4}x^4 y
\right)
\\
\textstyle
\kappa_1^{(3)}=\kappa_1^{(4)}=0
&\textstyle
\quad\Longrightarrow\quad
&
j^6 f(0)
\sim 
\left(x, y^2, xy^3-\frac{\kappa_1^{(5)}}{40\tau_0'\tau_1^5}x^5 y
\right)
\end{array}
$$
Since $C_k^\pm: (x, y^2, xy^3\pm x^ky)$ is $(k+1)$-determined, 
we see that 
$f$ is of type $C_k^\pm$ $(k=4,5)$
if and only if 
$\tau_0=\kappa_1=\kappa_1'=\cdots=\kappa_1^{(k-1)}=0$
and 
$\kappa_1^{(k)}\tau_0'\tau_1\lessgtr 0$ (seemingly, it is so for any $k$). 

$\bullet \;\; $ 
If $\kappa_1^{(3)}=2\tau_1\tau_0'\neq 0$, 
we have $j^4f(0)\sim (x, y^2, x^3 y)$ and 
$$f \sim_\A F_4: (x, y^2, x^3y+y^5) 
\Longleftrightarrow 
3\kappa_1^{(4)}-8 \tau_0' \tau_1'-12  \tau_1 \tau_0''\not=0.$$
All cases of sufficient jets of $\A$-simple germs have been covered. 
\end{rem}

\begin{rem}\label{5jet}\upshape{\bf (Non-realizable jets)} 
Let us continue the argument in Remark \ref{C_kF_4}. 
If  $\kappa_1^{(3)}=\tau_0'=0$, then $f$ should be of codimension $\ge 7$ and 
a computation shows that 

$$ 
j^5 f(0)\sim 
\left(x, y^2, \kappa_1^{(4)} x^4 y
+(\kappa_1^{(4)}-4\tau_1\tau_0'' )y^5
+2\sqrt{5}(\kappa_1^{(4)}-2\tau_1\tau_0'') x^2 y^3\right). 
$$ 
In particular, 
if two of 
three coefficients $\kappa_1^{(4)}$, $\kappa_1^{(4)}-4\tau_1\tau_0''$, 
$\kappa_1^{(4)}-2\tau_1\tau_0''$ are zero, 
then all are zero. 
Thus, for instance, the following $5$-jets are not equivalent to 
jets of any non-cylindrical ruled surface: 
$$
(x, y^2, x^4 y),
\quad(x, y^2, x^2y^3),
\quad(x, y^2, y^5).
$$
The $5$-jet $(x, y^2, y^5)$ is obviously realizable by a cylinder, 
while 
the $5$-jets $(x, y^2, x^4 y)$ and $(x, y^2, x^2y^3)$ are not equivalent to jets of 
{\it any} ruled surfaces, even if we drop the condition $\be'(0)\not=0$. 
In fact, put $F=\br(s)+t \be(s)$ with $\br(s)\cdot \be(s)=0$ and 
$\be(s)=(1, 0, 0) + o(s)$.  
If $F$ is singular at $(s,t)=(0,0)$ and $\br(0)=0$,  then $\br(s)=o(s)$. 
It is easy to see that 
$F\sim_\A f=(x, y^2h(x,y), y^3g(x,y))$ 
with some functions $h, g$ of the form $p(y)+xq(y)$, and 
thus the $5$-jet of $F$ is never equivalent to those two jets mentioned above. 
By the same reason, 
the $\A^3$-orbit of the $3$-jet $(x, y^3, x^2y)$ is not realized 
by jets of any ruled surfaces 
(the $2$-jet $(x,0,0)$ never appears in non-cylindrical ruled surfaces as seen before, and 
the $3$-jet is not realizable by ruled surfaces with $\be'(0)=0$, that is seen in the same way as above). 
\end{rem}

\subsection{Transversality} \label{mapping_space}
To precisely state {\it genericity} of  ruled surfaces, 
we need an appropriate mapping space (moduli space) equipped with a certain topology. 
By the definition, 
a {\it residual} subset of a mapping space is a union of countably many open dense subsets. 
When maps having a prescribed condition form a residual subset, 
we often say that such a map is {\it generic}, abusing words. 
Let $I$ be an open interval containing $0 \in \R$ and 
let $u$ denote the coordinate of $I$. 
As the mapping space of non-cylindrical ruled surfaces, we take 
$$\mathcal{R}:=\{\, \cbl=\bv_0+\varepsilon\bv_1 \in C^\infty(I, \cbU) \, | \, 
\bv_0'(u)\not=0\, (u \in I) \, \}$$ 
equipped with Whitney $C^\infty$ topology 
(As a remark,  Izumiya and Takeuchi \cite{IT} and Martins and Nu\~no-Ballesteros \cite{MN} 
took the space $C^\infty(I, \R^3\times S^2)$ instead of $C^\infty(I, \cbU)$,  
but the difference does not affect the matter of genericity arguments -- 
given a pair $(\br, \be)$ of base and director curves, we simply assign a curve $\cbl: I \to \cbU$ with $\bv_0=\be$ and $\bv_1=\br\times \be$). 

Also we put 
$$\mathcal{M}:=C^\infty(I, \R_{>0}\times \R^3)$$ 
of quadruples $(\kappa_0, \kappa_1, \tau_0, \tau_1)$ of real-valued functions with 
$\kappa_0(u)>0$ equipped with Whitney $C^\infty$ topology. 
Any curve $\cbl(u)$ in  $\mathcal{R}$ defines $\bD$-valued functions, 
$\check{\kappa}(u)$ and  $\check{\tau}(u)$ (parameterized by a general parameter $u \in I$), 
that produces a continuous map $\Phi: \mathcal{R} \to \mathcal{M}$. 
Obviously, $\Phi$ is surjective. In fact, given 
a quadruple of functions  $(\kappa_0(u), \kappa_1(u), \tau_0(u), \tau_1(u)) \in \mathcal{M}$, 
put a new parameter $s:=s(u)=\int_0^u \kappa_0(u) du$ 
and define $\kappa_1(s):=\kappa_1(u(s))$, etc. 
Then, three functions $\kappa_1(s), \tau_0(s), \tau_1(s)$ 
determines, up to Euclidean motions,  the curve $\cbl(s)=\bv_0(s)+\varepsilon \bv_1(s)$ 
by solving the ordinary differential equation 
determined by the Frenet formula.  
The ambiguity is  
fixed by the initial values $\cbl(0), \cbn(0), \cbt(0)$, 
which corresponds to the initial orthogonal axes in $\R^3$ 
at $u=0$. 
Put $\cbl(u):=\cbl(s(u)) \in \mathcal{R}$; the set of such cruves 
is exactly the preimage via $\Phi$ of the given quadruple of functions. 
That infers that for a dense subset $O \subset \mathcal{M}$, 
the preimage $\Phi^{-1}(O)$ is also dense in $\mathcal{R}$. 

The above construction is extended for a parametric version. 
Let $W$ be an open subset of $\R^p\; (0\le p \le 3)$, 
and  consider the subspace 
$\mathcal{R}_W$ of $C^\infty (I\times W, \cbU)$ 
which consists of maps $\cbl(u, \lambda)=\bv_0(u, \lambda)+\varepsilon \bv_1(u, \lambda)$  
with parameter $\lambda \in W$ 
satisfying $\rd \bv_0/\rd u\not=0$ at any $(u, \lambda)$. 
Put $\mathcal{M}_W$ to be 
the mapping space of $I \times W \to \R_{>0}\times \R^3$, 
and then 
a surjective continuous map $\Phi:\mathcal{R}_W \to \mathcal{M}_W$ 
is defined in entirely the same way as above. 
For a dense subset $O \subset \mathcal{M}_W$, 
the preimage $\Phi^{-1}(O)$ is also in $\mathcal{R}_W$. 

As seen in the previous section, we have obtained 
a semi-algebraic stratification of the jet space $J^r:=\R^3 \times J^r(1,3)$ 
up to codimension $4$ ($r$ sufficiently large). 
In fact, 
any strata are defined by the conditions in Table \ref{table1} 
of (in)equalities in Taylor coefficients $\{\kappa_1^{(k)}, \tau_0^{(k)}, \tau_1^{(k)}\}_{0 \le k\le r}$, 
which form a system of coordinates of the affine space $J^r$. 
Notice that 
these Taylor coefficients are with respect to the arclength parameter $s$. 
For each quadruple $(\kappa_0, \kappa_1, \tau_0, \tau_1) \in \mathcal{M}_W$, 
we put  
$$s=s(u, \lambda):=\int_0^u \kappa_0(u, \lambda) du, \qquad  
\varphi(u, \lambda)=(\kappa_1(u,\lambda), \tau_0(u,\lambda), \tau_1(u,\lambda)).$$ 
By the assumption that $\rd s/\rd u=\kappa_0>0$, 
let $\bar{\varphi}(s, \lambda):=\varphi(u(s, \lambda), \lambda)$. 
Then we define 
$$\Psi: I \times W \times \mathcal{M}_W \to J^r, 
\qquad \Psi(u, \lambda, (\kappa_0, \varphi)) := 
j^r_s\bar{\varphi}(s(u, \lambda), \lambda),$$
where $j^r_s\bar{\varphi}$ means the $r$-jet respect to the parameter $s$. 
By a version of Thom's transversality theorem (Lemma 4.6 in \cite{GG}), 
there is a dense subset $O$ of $\mathcal{M}_W$ 
so that for any $\varphi \in O$, 
the jet extension $\Psi_{\kappa_0, \varphi}: I \times W \to J^r$ 
is transverse to every stratum of our stratification of $J^r$. 
Hence, $\Phi^{-1}(O)$ is dense in $\mathcal{R}_W$, 
and for any element of $\Phi^{-1}(O)$, 
only $\A$-singularity types listed in Table \ref{table1} appears. 
This completes the proof of (2) in Theorem \ref{thm1}.  
\qed 

\begin{rem}\label{versal}\upshape 
{\bf ($\A_e$-versal deformations)}
For each type in Table \ref{table1}, 
an $\A_e$-versal deformation of the germ is realized 
by a generic family of non-cylindrical ruled surfaces. 
This is directly checked by computations. 
For instance, as in Table \ref{table1}, 
the $S_1^\pm$-singularity of ruled surface at $s=0$ is characterized by 
$\kappa_1(0)=\kappa_1'(0)=0$,  
$\kappa_1''(0)\not=0, 2\tau_0(0)\tau_1(0)$ and $\tau_1(0)\not=0$. 
Suppose that $\varphi=(\kappa_1(s), \tau_0(s), \tau_1(s)): I \to \R^3$ satisfies this condition. 
Define a $1$-parameter family $I \times \R \to \R^3$ 
by $\varphi(s,\lambda) :=\varphi(s)+(\lambda, 0, 0)$, then 
obviously, its $1$-jet extension 
$j^1_s\varphi$ is transverse at $(0,0)$ to 
the stratum defined by $\kappa_1=\kappa_1'=0$ in $J^1=\R^3\times J^1(1,3)$. 
This family yields a $1$-parameter family 
$F(s, t, \lambda)=(t, ts - \frac{\tau_1}{2}s^2, \lambda s)+o(2)$ 
of ruled surfaces. 
By using a coordinate change of $x=t+\cdots$ 
(= first component of $F(s, t, \lambda)$) and $y=s$ and some target changes, 
we see that 
the germ of $F(s, t, \lambda)$ is equivalent to $(x, y^2, y^3\pm x^2y+\lambda y)$, 
which is an $\A_e$-miniversal deformation of $S_1^\pm$-singularity. 
\end{rem}

\section{Singularities of developable surfaces}

\subsection{Recognition of singularity types}
For non-cylindrical developable surfaces, $\kappa_1(s)\equiv 0$ identically. Hence 
the Taylor expansion of $f$ is (\ref{2}) with $\kappa_1^{(k)}= 0$ for all $k$: 
$$\textstyle
f(x,y):=F(y,t(x,y))=
(x, x y - \frac{1}{2}\tau_1 y^2 - \frac{1}{6}\tau_1' y^3, 
\frac{1}{2}\tau_0 x y^2 + \frac{1}{3} \tau_0\tau_1 y^3)+o(3).$$
Using the $\A$-criteria mentioned in \S 2, 
we classify singularities arising in generic families of developable surfaces. 
Notice that there are two different aspects; 
 singularities of frontal surfaces correspond to the case of $\tau_1\not=0$, 
 while  singularities of wavefronts correspond to the case of $\tau_1=0$. 
 Below we prove Theorem \ref{thm2}. 

\

\t
$\bullet\;\; $ 
{\bf Case of $\tau_1\not=0$}: 
By $s = y + \tau_1^{-1}x$ and some linear change of the target, 
we have $f=(x, y^2+o(2), f_3(x,y))$ with $f_3=\tau_0y^3+o(3)$. 
Note that $(x, y^2)$ is $2$-determined and that 
each term $x^ky^{2l}$ in $f_3$ can be removed by a coordinate change of the target 
$(X, Y, Z) \mapsto (X, Y, Z-X^kY^l)$. 
Use Proposition \ref{CE} in \S 2 (\cite{Mond89}) 
for determinacy in CE. 

\begin{itemize}
\item[(i)] 
If $\tau_0\not=0$,  then $f\sim_\A(x,y^2,y^3)$, since it is $3$-determined in CE. 
\item[(ii)] 
Let $\tau_0=0$. 
Computing the $4$-jet,  we see 
$$f_3=\tau_0'(6x^2y^2+8\tau_1xy^3+3\tau_1^2y^4)+o(4).$$  
If $\tau_0'\not=0$, then 
$f\sim_\A (x, y^2,  x y^3)$, for 
the germ is $4$-determined in CE. 
Hence $f$ is of type cuspidal crosscap. 
\item[(iii)] 
Let $\tau_0=\tau_0'=0$. 
Computing the $5$-jet,  we see 
$$f_3=
\tau_0''(10x^3y^2+20\tau_1x^2y^3+15\tau_1^2xy^4+4\tau_1^3y^5)+o(5).$$  
If $\tau_0''\not=0$, by target changes using $X=x$ and $Y=y^2$, 
terms $x^3y^2$ and $xy^4$ can be removed from $Z=f_3$, 
thus we see that 
$f\sim_\A (x, y^2, y^3(x^2+y^2))$, 
for this germ is $5$-determined in CE. 
That is cuspidal $S_1^+$-type. Note that cuspidal $S_1^-$ never appears. 
\item[(iv)] 
Let $\tau_0=\tau_0'=\tau_0''=0$.  
Computing the $6$-jet,  we see 
$$f_3=\tau_0'''(15x^4y^2+40\tau_1x^3y^3
+45\tau_1^2x^2y^4+24\tau_1^3xy^5+5\tau_1^4y^6)+o(6).$$  
If $\tau_0'''\not=0$, then 
$f\sim_\A (x, y^2, y^3(x^3+xy^2))$, 
for the germ is $6$-determined in CE. 
That is cuspidal $C_3^+$-type, while cuspidal $C_3^-$ does not appear. 
Note that $\tau_0=\tau_0'=\tau_0''=0$ if and only if the $5$-jet of $f$ is equivalent to $(x, y^2,0)$, 
thus cuspidal $S$ and $B$-types never appear, as mentioned in Remark \ref{rem1}. 
\end{itemize}

\

\t
$\bullet\;\; $
{\bf Case of $\tau_1=0$}:  
Then 
$\textstyle 
f=(x, x y - \frac{1}{6}\tau_1'y^3, \frac{1}{2}\tau_0 x y^2) + o(3)$. 
Note that $j^2f(0) \sim (x, xy,0)$, thus 
types $A_3^\pm$ and $D_k$ never appear (Remark \ref{rem1}). 

If $\tau_0=0$,  
$j^3f(0)$ is equivalent to either $(x, xy+y^3,0)$ or $(x, xy,0)$, 
that is of type $T_1$ or $T_2$ (codimension $3, 4$) in Table \ref{table2}. 
Now assume that $\tau_0\not=0$. 
Write 
$$\textstyle
f=(x, f_2(x,y), f_3(x,y))=(x, x y - \frac{1}{6}\tau_1'y^3,  x y^2)+o(3).$$ 
The singular point set $S(F)$ is defined by $(f_2)_y=(f_3)_y=0$, 
and through a computation, it is simplified as $\lambda=0$ with 
$$\textstyle 
\lambda=x-\frac{1}{2}\tau_1'y^2 - \frac{1}{6}\tau_1''y^3
-\frac{1}{24}(\tau_1'''-3\tau_1')y^4+o(4).$$
We may take $\eta={\rd}/{\rd y}$ as a vector field which generates $\ker dF$ along $S(F)$. 
Then,  $\eta \lambda(0)=0$,   
$\eta\eta \lambda(0) = -\tau_1'$, 
$\eta\eta\eta \lambda(0) =- \tau_1''$ and 
$\eta\eta\eta\eta \lambda(0) =-(\tau_1'''-3\tau_1')$. 
Hence, by Izumiya-Saji's criteria in \S 2.5, 
we have the conditions for detecting $Sw$, $cA_4$ and $cA_5$.

\subsection{Topological classification}
We prove Theorem \ref{thm3}. 
Let $\sigma(s)$ be the striction curve of a non-cylindrical developable surface. 
Assume that $\sigma(0)=0 \in \R^3$, and consider the germ $\sigma:(\R,0)\to (\R^3,0)$. 
Since  $\{\bv_0(s), \bn_0(s), \bt_0(s)\}$ form a basis of $\R^3$ 
for each $s$,  we denote  the $k$-th derivative by 
$$
\sigma^{(k)}(s)
=
A_k(s)\bv_0(s)
+B_k(s)\bn_0(s)
+C_k(s)\bt_0(s)
\quad(k\geq 1)
$$
where $A_k(s), B_k(s), C_k(s)$ are some functions. 
Then, with respect to the basis $\{\bv_0(0), \bn_0(0), \bt_0(0)\}$, 
the expansion of $\sigma$ at $s=0$ 
is given by 
$$\textstyle 
\sigma(s)=(A_1(0)s+\frac{1}{2}A_2(0)s^2+\cdots, 
B_1(0)s+\frac{1}{2}B_2(0)s^2+\cdots, 
C_1(0)s+\frac{1}{2}C_2(0)s^2+\cdots).$$ 
Now assume that $\sigma$ is of type $(m, n_1, n_2)$, i.e., 
$$
\begin{cases}
A_1(0)=\cdots=A_{m-1}(0)=0,\ A_m(0)\neq 0,
&\\
B_1(0)=\cdots=B_{n_1-1}(0)=0,\ B_{n_1}(0)\neq 0,
&\\
C_1(0)=\cdots=C_{n_2-1}(0)=0,\ C_{n_2}(0)\neq 0. 
&\end{cases}
$$
Since $\sigma'(s)=\tau_1(s)\bv_0(s)$ for a developable surface (Lemma \ref{striction_curve0} (iii)), 
we see that $A_1(s)=\tau_1(s)$ and $B_1(s)\equiv C_1(s)\equiv 0$. 
By the Frenet formula (Theorem \ref{frenet} (1)), 
\begin{align*}
\sigma^{(k+1)}
&=(\sigma^{(k)})'
=\{A_{k}\bv_0
+B_{k}\bn_0
+C_{k}\bt_0\}'\\
&=
(A_{k}'-B_k)\bv_0
+(B_{k}'+A_k- C_k \tau_0)\bn_0
+(C_{k}' +B_k \tau_0)\bt_0 \\
&=A_{k+1}\bv_0+B_{k+1} \bn_0 + C_{k+1} \bt_0. 
\end{align*}
Thus for $k=1$, we have 
$A_2(s)=\tau_1'(s)$, $B_2(s)=\tau_1(s)$, $C_2(s)\equiv 0$, and 
for $k=2$, 
$A_3(s)=\tau_1''(s)-\tau_1(s)$, 
$B_3(s)=2\tau_1'(s)$ and $C_3(s)=\tau_0(s)\tau_1(s)$. 
For $k \ge 3$, 
there are some smooth functions 
$a_{k,*}(s)$, $b_{k,*}(s)$, $c_{k,*,*}(s)$ and 
positive numbers 
$\beta_k, \gamma_{k,0}, \cdots, \gamma_{k, k-3}>0$ such that 
\begin{align*}
A_{k}(s)
&=
a_{k,0}(s)\tau_1(s)
+\cdots
+a_{k,k-2}(s)\tau_1^{(k-2)}(s)
+
\tau_1^{(k-1)}(s),
\\
B_{k}(s)
&=
b_{k,0}(s)\tau_1(s)
+\cdots
+b_{k,k-3}(s)\tau_1^{(k-3)}(s)
+\beta_k\tau_1^{(k-2)}(s),
\\
C_{k}(s)
&=
\{c_{k,0,0}(s)\tau_0(s)+\cdots
+\gamma_{k,0}\tau_0^{(k-4)}(s)\}\tau_1(s)
\\
&\quad
+\{c_{k,1,0}(s)\tau_0(s)+\cdots
+\gamma_{k,1}\tau_0^{(k-5)}(s)\}
\tau_1'(s)
+\cdots
\\
&\quad
+\{c_{k,k-4,0}(s)\tau_0(s)+\gamma_{k,k-4}\tau_0'(s)\}
\tau_1^{(k-4)}(s)
+\gamma_{k,k-3}\tau_0(s)
\tau_1^{(k-3)}(s). 
\end{align*}
Hence, by the assumption on $A_k(0)$, we have 
$$\tau_1(0)=\cdots=\tau_1^{(m-2)}(0)=0,\quad
\tau_1^{(m-1)}(0)\neq 0,$$
and thus 
$$B_1(0)=\cdots=B_m(0)=0, \;\; B_{m+1}(0)\neq 0, 
\quad C_1(0)=\cdots=C_{m+2}(0)=0.$$
In particular, 
$$n_1=m+1, \;\; n_2=m+1+r\quad (r\geq 1).$$
By the above formula of $C_k(s)$ with $k=m+1+r$, 
we see 
$$
\tau_0(0)=\cdots=\tau_0^{(r-2)}(0)=0,
\quad 
\tau_0^{(r-1)}(0)\neq 0. 
$$
Conversely, 
if the order of $\tau_0$ and $\tau_1$ 
are $r$ and $m-1$, respectively, 
then the type of $\sigma$ is $(m, m+1, m+1+r)$. 
This completes the proof. \qed



\begin{thebibliography}{99999}
%
%
\bibitem{GG} M. Golubitsky and V. Guillemin, 
\textit{Stable Mappings and Their Singularities}, 
GTM 14, Springer-Verlag (1973). 
%
\bibitem{Guggenheimmer} H. W. Guggenheimer, 
\textit{Differential Geometry}, 
McGraw-Hill (1963). 
%
\bibitem{Hlavaty} V. Hlavar\'y, 
\textit{Differential line geometry}, Noordhoff, Groningen-London (1953). 
%
\bibitem{HT} S. Honda and M. Takahashi, 
Framed curves in the Euclidean space, Adv. Geom. {\bf 16} (3) (2016), 265--276. 
%
\bibitem{Ishikawa}  G. Ishikawa, 
Topological classification of the tangent developables of space curves, 
J. London Math. Soc. (2) {\bf 62} (1999), 
583--598. 
%
\bibitem{Ishikawa99}  G. Ishikawa, 
Singularities of developable surfaces, 
Singularity Theory, Proc. European Sing. Conf. (Liverpool, 1996), 
ed. W. Bruce and D. Mond, Cambridge Univ. Press (1999), 
403--418. 
%
\bibitem{IS} S. Izumiya and K. Saji, 
The mandala of Legendrian dualities for pseudo-spheres in 
Lorentz-Minkowski space and ``flat" spacelike surfaces, 
J. Singularities {\bf 2} (2010), 92--127. 
%
\bibitem{IST} S. Izumiya, K. Saji and M. Takahashi, 
Horospherical flat surfaces in Hyperbolic $3$-space, 
J. Math. Soc. Japan {\bf 62} (2010), 789--849. 
%
\bibitem{IT} S. Izumiya and N. Takeuchi, 
Singularities of ruled surfaces in $\R^3$, 
Math. Proc. Camb. Phil. Soc. {\bf130} (2001), 1--11. 
%
\bibitem{ICRT} S. Izumiya, M. C. Romero Fuster, M. A. S. Ruas and F. Tari, 
\textit{Differential Geometry from a Singularity Theory Viewpoint}, 
 World Scientific (2016). 
%
\bibitem{Kabata}  Y.~Kabata,
Recognition of plane-to-plane map-germs, 
Topology and its Appl. {\bf 202} (2016), 216--238. 
%
\bibitem{Kurokawa}  H. Kurokawa, 
On generic singularities of $1$-parameter families of developable surfaces  (in japanese), 
Master Thesis, Hokkaido University (2013). 
%
\bibitem{MN} R. Martins and J. J. Nu\~no-Ballesteros, 
Finitely determined singularities of ruled surfaces in $\R^3$, 
Math. Proc. Camb. Phil. Soc. {\bf 147} (2009), 701--733. 
%
\bibitem{Mond82} D. Mond, 
\textit{Classification of certain singularities and applications to differential geometry}, 
Ph.D. thesis, The University of Liverpool (1982).
%
\bibitem{Mond85} D. Mond, 
On the classification of germs of maps from $\R^2$ to $\R^3$, 
Proc. London Math. Soc. (3) {\bf 50} (1985), 333--369. 
%
\bibitem{Mond89} D. Mond, 
Singularities of the tangent developable surface of a space curve, 
Quart. J. Math. Oxford Ser. (2), {\bf 40} (1989), 79--91. 
%
\bibitem{PW}  H. Pottmann and J.  Wallner, 
\textit{Computational Line Geometry},   
Mathematics and Visualization, Springer (2001). 
%
\bibitem{Saji} 
K. Saji, Criteria for $S_k$ singularities and their applications, 
Jour. G\"okova Geom. Top. {\bf 4} (2010), 67--81. 
%
\bibitem{Salmon} G. Salmon,
\textit{A treatise on the analytic geometry of three dimensions},
4th edition, Dublin (1882).
%
\bibitem{SKSO} H. Sano, Y. Kabata, J. L. Deolindo Silva and T. Omoto, 
Classification of jets of surfaces in $3$-space via central projection, 
Bull. Brasilian Math. Soc. New Series (2017). 
DOI 10.1007/s00574-017-0036-x. 
%
\bibitem{Selig} J. M. Selig, 
\textit{Geometric Fundamentals of Robotics, 2nd edition}, 
Monographs in Computer Science, Springer (2005). 
%
\bibitem{Tanaka}  J. Tanaka, 
Clifford algebra and singularities of ruled surfaces (in Japanese), 
Master Thesis, Hokkaido University (2016). 
\end{thebibliography}
\end{document}